\documentclass[EJP]{ejpecp}

\pdfoutput=1

\AUTHORS{%
  Magdalena~Kobylanski\footnote{Universit\'e Paris-Est, UMR 8050 (LAMA),
    France.~\EMAIL{magdalena.kobylanski@univ-mlv.fr}
}
  \and 
  Marie-Claire~Quenez\footnote{Universit\'e Denis Diderot (P7), UMR 7599 (LPMA), INRIA, France.
\EMAIL{quenez@math.jussieu.fr}
}}

\SHORTTITLE{Optimal stopping in a general framework}

\TITLE{Optimal stopping time problem in a general framework}

\KEYWORDS{optimal stopping ; supermartingale ; american options}

\AMSSUBJ{60G40} 

\SUBMITTED{August 29, 2011}
\ACCEPTED{July 27, 2012} 
\VOLUME{17}
\YEAR{2012}
\PAPERNUM{72}
\DOI{I10.1214/EJP.v17-2262}
\ABSTRACT{We study the optimal stopping time problem $v(S)={\rm ess}\sup_{\theta \geq S} E[\phi(\theta)|\mathcal
{F}_S]$, for any stopping time $S$,  where the reward is given by a family $(\phi(\theta),\theta\in\mathcal{T}_0)$ 
\emph{of non negative random variables} indexed by stopping times. We solve the problem under weak assumptions in terms of integrability and regularity of the reward family. More precisely, we only suppose $v(0) < + \infty$ and $
(\phi(\theta),\theta\in \mathcal{T}_0)$ upper semicontinuous along stopping times in expectation. We show the 
existence of an optimal stopping time and obtain a characterization of the minimal and the maximal optimal stopping 
times. We also provide some local properties of the value function family. All the results are written in terms of families of random variables and are proven by only using classical results of the Probability Theory. }
\begin{document}
\section*{Introduction}

In the present work we study the optimal stopping problem in the setup of families of random variables indexed by  stopping times, which is more general than the classical setup of processes. This allows technically simpler and clearer proofs, and also to solve the problem under weaker assumptions. 

To the best of our knowledge, the most general result given in the literature is that of El Karoui (1981):  existence of an optimal stopping time is proven when the reward is given by an upper semicontinuous non negative process of 
class $\mathcal{D}$.  For a classical exposition of the Optimal Stopping Theory, we also refer to Karatzas Shreve (1998) and Peskir Shiryaev (2005), among others.

\vline

Let $T\in \mathbb{R}_+^*$ be the terminal time and let $(\Omega,\mathcal{F},  (\mathcal{F}_t)_{0\le t\le T}, P)$ be a filtered probability set which satisfies the usual  conditions.

 An optimal stopping problem can be naturally expressed in terms of families of random variables  indexed by stopping times. Indeed,
consider an agent who can choose a stopping time in $\mathcal{T}_0$. When she decides to stop at $\theta$ 
$\in$ $\mathcal{T}_0$, she receives the amount $\phi(\theta)$, where $\phi(\theta)$ is a non negative $\mathcal{F}_\theta$-measurable random variable. The family  $(\phi
(\theta), \theta\in\mathcal{T}_0)$ of random variables indexed by stopping times is called the \emph{reward (or payoff) family}. It is identified with the map $\phi:\theta \mapsto  \phi(\theta)$   from $\mathcal{T}_0$ into the set of random variables.

A family $(\phi(\theta), \theta\in\mathcal{T}_0)$   is said to be  an \emph{admissible} family if it satisfies the two following conditions. First,  for each stopping time $\theta$,   $\phi(\theta)$ is  a non negative and  $\mathcal{F}_\theta$-measurable random variable. Second, the following natural compatibility condition holds: for each $\theta$, $\theta'$ in  $\mathcal{T}_0$, $\phi(\theta)=\phi(\theta')$ a.s. on the subset $\{\theta=\theta'\}$ of $\Omega$.

 In the sequel, the \emph{reward  family}  $(\phi(\theta), \theta\in\mathcal{T}_0)$ is supposed to be an \emph{admissible} family of non negative random variables.

At time $0$, the agent wants to choose a stopping time $\theta^*$ so that it maximizes $E[\phi(\theta)]$ over the set of stopping times $\mathcal{T}_0$. 
The best expected reward at time $0$ is thus given by $v(0):=\sup_{\theta \in\mathcal{T}_0} E[\phi(\theta)]$,  and is also called  the value function at time $0$. Similarly, for  a stopping time $S$ $ \in\mathcal{T}_0$, the \emph{value function at  time $S$} is defined by 
\[v(S):={\rm ess} \sup\{\, E[\phi(\theta)\, |\,\mathcal{F}_S], \;\theta \in\mathcal{T}_0\;{\rm and}\; \theta \geq S\,{\rm a.s.}\}.\]
The  family of random variables $v=(v(S), S\in\mathcal{T}_0)$ can be shown to be admissible, and characterized as the \emph{Snell envelope family} of $\phi$, also denoted by $\mathcal{R} (\phi)$,  defined here as the smallest supermartingale family greater than the reward family $\phi$.
The Snell envelope operator $\mathcal{R}: \phi\mapsto\mathcal{R} (\phi)=v
$, thus acts on the set of admissible families of r.v. indexed by stopping times.

  Solving the optimal stopping time problem at time $S$ mainly consists to prove the existence of an optimal stopping times $\theta^*(S)$, that is, such that $v(S)= E[\phi(\theta^*(S))|\mathcal{F}_S]$ a.s.\,
  
 Note that this setup of families of random variables indexed by stopping times is  clearly more general than the setup of processes. Indeed, if $(\phi_t)_{0\le t\le T}$ is a progressive process,  we set  $\overline \phi (\theta):=\phi_\theta$, for each stopping time $\theta$. Then,   the family $\overline\phi =(\phi_\theta,\theta\in\mathcal{T}_0)$ is admissible.  

\vline

The interest of such families has already been stressed, for instance, in the first chapter of El Karoui (1981).
However, in that work as well as in the classical litterature, the optimal stopping time problem is set and solved in the setup of processes. In this case, the reward is given by a progressive process $(\phi_t)$ and  the associated value function family $(v(S), S\in\mathcal{T}_0)$ is defined as above but does not \emph{a priori} correspond to a progressive process. 
An important  step of the classical approach consists in \emph{aggregating} this familly, that is in finding a process $(v_t)_{0\le t\le T}$ such that, for each stopping time $S$, $v(S)=v_S$ a.s.\, This aggregation problem is solved by using some fine results of the General Theory of Processes.
Now, it is well known that this process $(v_t)$ is also  characterized as the \emph{Snell envelope process} of the reward process $(\phi_t)$. Consequently, the previous aggregation result allows to define the \emph{Snell envelope process operator} $\hat{\mathcal{R}}:(\phi_t)\mapsto \hat{\mathcal{R}} [(\phi_t)]: = (v_t)$ which acts here on the set of progressive processes. The second step then consists, by a  penalization method introduced by Maingueneau (1978), and under some right regularity conditions on 
the reward process, in showing the existence of $\varepsilon$-optimal stopping time. Next, under additional left regularity conditions on the 
reward process,  the minimal optimal stopping time is characterized as a hitting time of processes, namely
 \begin{equation*}\label{optun}
\overline \theta (S):=\inf\{\,t\geq S, \;v_t=\phi_t\,\}\,.
 \end{equation*}
 Finally, as the value function $(v_t)$ is a strong supermartingale of class $\mathcal{D}$, it admits a \emph{Mertens decomposition}, which, in the right continuous case, reduces to the \emph{Doob-Meyer decomposition}. This decomposition is then used to characterize the maximal optimal stopping time as well as to obtain some local properties of the value function $(v_t)$. The proofs of these properties thus rely on strong and sophisticated results of the General Theory of Processes (see the second chapter of El Karoui (1981) for details). It is not the case in the framework of admissible families.

\vline

 In the present work, which is self-contained, we study the general case of a reward given by an \emph{admissible family}  $\phi=(\phi(\theta), \theta\in\mathcal{T}_0)$
 of \emph{non negative} random variables, and we solve the associated optimal stopping time problem only in terms of admissible families. Using this approach, we avoid the \emph{aggregation step} as well as the use of \emph{Mertens' decomposition}.  
 
 Moreover, we only make the assumption $v(0)= \sup_{\theta\in\mathcal{T}_0} E[\phi(\theta)] <+ \infty$, which is, in the case of a reward process, weaker than the assumption $(\phi_t)$ of class $\mathcal{D}$, required in the previous literature.

Furthermore,  the existence $\varepsilon$-optimal stopping times is obtained when  $(\phi(\theta), \theta\in\mathcal{T}_0)$ is right upper semicontinuous along stopping times in expectation, that is, for each stopping time $\theta$, and, for each non decreasing sequence $(\theta_n)$ of stopping time tending to $\theta$,  $\limsup_n E[\phi(\theta_n)]\le E[\phi(\theta)]$. This condition is, in the case of a reward process, a bit wilder than the usual assumption ``$(\phi_t)$ right upper semicontinuous and of class $\mathcal{D}$''.

Then, under the additional assumption that the reward family is left upper semicontinuous along stopping times in expectation, we show the existence of  optimal stopping times  and we characterize the minimal optimal stopping time $\theta_{*}(S)$ for $v(S)$  by 
 \[\theta_{*}(S) ={\rm ess} \inf \{\,\theta \in\mathcal{T}_0,\,\, \theta \geq S\;{\rm a.s.} \,{\rm and} \; u(\theta) =\phi(\theta) \,\,\mbox {\rm a.s.} \,\}.\]
Let us emphasize that $\theta_*(S)$ is no longer defined as a hitting time of processes but as an essential infimum of a set of stopping times. This formulation is a key tool to solve the optimal stopping time problem in the unified framework of admissible families.

Furthermore, we introduce the following random variable \begin{equation*}
\check{\theta}(S) :={\rm ess} \sup\{\;\theta \in\mathcal{T}_0,\,\, \theta \geq S\,\,{\rm a.s.}\,\,{\rm and}\,\, E[v(\theta)]= E[v(S)]  \,\},
\end{equation*}
and show that it is the maximal optimal stopping time for $v(S)$. 

Some local properties of the value function family $v$ are also investigated. To that purpose, some new local notions for families of random variables are introduced.
We point out that  these properties are proved using only classical probability results. 
In the case of processes, these properties correspond to some known results shown, using very sophisticated tools, by Dellacherie and Meyer (1980) and El Karoui (1981), among others.

 \vline
 
 At last, let us underline that the setup of families of random variables indexed by stopping time was used by Kobylanski et al. (2011), in 
 order to study  optimal multiple stopping. This setup is particularly relevant in that case. In particular, it avoids the 
 aggregation problems, which, in the case of multiple stopping times, appear to be particularly knotty and difficult.  The setup of families of random variables  is also used in Kobylanski et al. 
 (2012) to revisit the Dynkin game problem and provides a new insight on this well-known problem. 

\vline

Let $\mathbb{F}=(\Omega,\mathcal{F}, (\mathcal{F}_t)_{0\le t\le T},P)$ be a probability space which
filtration 
$(\mathcal{F}_t)_{0\le t\le T}$ satisfies the usual conditions of right continuity and augmentation by the null sets of $\mathcal{F}=
\mathcal{F}_T$. We suppose that $\mathcal{F}_0$ contains only sets of probability $0$ or $1$.  The time horizon is a fixed constant  
$T\in ]0,\infty[$.
We denote by $\mathcal{T}_{0}$ the collection of stopping times of
${\mathbb{F} }$ with values in $[0 , T]$.  More generally, for any
stopping times $S$, we denote by $\mathcal{T}_{S}$  (resp.  $\mathcal{T}_{S^+}$) the class of stopping times
$\theta\in\mathcal{T}_0$ with $\theta\geq S$ a.s.\, (resp. $\theta>S $ a.s. on $\{S<T\}$ and $\theta=T$ a.s. on $\{S=T\}$).

For $S,S'$ $\in $ $\mathcal{T}_0$, we also define $\mathcal{T}_{[S, S']}$ the set of $\theta\in\mathcal{T}_0$ with $S \leq \theta \leq S'$ a.s. and $\mathcal{T}_{]S, 
S']}$  the set of $\theta\in\mathcal{T}_0$ with $S < \theta \leq S'$ a.s.

Similarly, the set ``$\mathcal{T}_{]S, S']}$ on $A$''  denotes the set of $\theta\in\mathcal{T}_0$ with $S < \theta \leq S'$ a.s.\, on $A$.

We use the following notation: for real valued  random variables $X$ and $X_n$, $n\in$ $\mathbb N$,  ``$X_n\uparrow X$'' stands for ``the sequence $(X_n)$ is nondecreasing and converges to $X$ a.s.''.  
\section{First properties}
In this section we prove some results about the value function families $v$ and $v^+$ when the reward is given by an admissible  family of random variables indexed by stopping times.   Most of these results are, of course, well-known in the case of processes. 
\begin{definition}\label{def.admi}
We say that a family 
$\phi=(\phi(\theta), \, \theta\in \mathcal{T}_0)$  is \emph{admissible} if it satisfies the following conditions 
\par
1. \quad for all
$\theta\in \mathcal{T}_0$ $\phi(\theta)$ is a $\mathcal{F}_\theta$-measurable non negative random variable,
\par
 2. \quad  for all
$\theta,\theta'\in \mathcal{T}_0$, $\phi(\theta)=\phi(\theta')$ a.s.  on
$\{\theta=\theta'\}$.
\end{definition}
\begin{remark}
By convention, the non negativity property of a random variable means that it takes its values in $\overline {\mathbb{R}} ^+$.

Also, it is always possible to define a admissible family associated with a given process.
More precisely, let $(\phi_t)$ be a non negative progressive process. Set $\overline \phi(\theta):= \phi_\theta$, for each $\theta$ $\in$ $\mathcal{T}_0$ . Then, the family $\overline\phi= (\phi_\theta, \, \theta\in \mathcal{T}_0)$ is  clearly admissible.
\end{remark}
 Let $(\phi(\theta), \, \theta\in \mathcal{T}_0)$ be an admissible family called \emph{reward}. For $S\in \mathcal{T}_0$, the \emph{value function at time $S$} is defined by
\begin{equation}\label{eq.vs}
v(S):={\rm ess} \sup_{\theta\in \mathcal{T}_S} E[\phi(\theta) \, |\,\mathcal{F}_S]\, ,
\end{equation}
the \emph{strict value function at time $S$} is defined by
\begin{equation}\label{eq.vsp}
v^+(S):={\rm ess} \sup_{\theta\in \mathcal{T}_{S^+}} E[\phi(\theta) \, |\,\mathcal{F}_S] \,.
\end{equation}
where $\mathcal{T}_{S^+}$ is the class of stopping times
$\theta\in \mathcal{T}_0$ with $\theta>S $ a.s. on $\{S<T\}$ and $\theta=T$ a.s. on $\{S=T\}$.
Note that $v^+(S)= \phi (T)$ a.s. on $\{S=T\}$.

\vline

Note that the essential supremum  of a family $\mathcal X$ of non negative random variables, denoted ``${\rm ess} \sup\,  \mathcal{X}$'',  is a well defined, almost surely unique random variable. Moreover, if  $\cal X$ is stable by pairwise maximization (that is $X\vee X'\in \mathcal{X}$ for all $X$ and $X'$ $\in$ $\mathcal{X}$), then there exists a sequence $(X_n)$ in $\cal X$ such that $X_n\uparrow ({\rm ess} \sup\,  \mathcal{X})$. We refer to Neveu (1975) for a complete and simple proof (Proposition VI-1.1. p 121).

\begin{proposition}\label{P1.Adm}\emph{(Admissibility of $v$ and $v^+$)}\\
The families $v=(v(S), S\in \mathcal{T}_0)$ and $v^+=(v^+(S), S\in \mathcal{T}_0)$ defined by (\ref{eq.vs}) and (\ref{eq.vsp}) are admissible.
\end{proposition}

\begin{proof} The arguments are the same for $(v(S), S\in \mathcal{T}_0)$ and  $(v^+(S), S\in \mathcal{T}_0)$.
We prove the property only for  $(v^+(S), S\in \mathcal{T}_0)$. 
Property 1 of admissibility for  $(v^+(S), S\in \mathcal{T}_0)$ follows from the existence of the essential supremum (see Neveu (1975)). \\
Take $S,S'\in \mathcal{T}_0$ and let $A=\{S=S'\}$. For each $\theta\in \mathcal{T}_{S^+}$ put $\theta_A=\theta \mathbb{1}_A+T \mathbb{1}_{A^c}$. As $A \in\mathcal{F}_S\cap\mathcal{F}_{S'}$, one has a.s. on $A$, 
$E[\phi(\theta)\,|\,\mathcal{F}_S]= E[\phi(\theta_A)\,|\,\mathcal{F}_S]=$ $E[\phi(\theta_A)\,|\,\mathcal{F}_{S'}]\leq$ $ v^+(S'),$ because $\theta_A$ $\in$ $\mathcal{T}_{S'^+}$.
Hence, taking the essential supremum over $\theta\in \mathcal{T}_{S^+}$ one has $v^+(S)\leq v^+(S')$ a.s. and by symmetry of $S$ and $S'$, we have proven property 2 of admissibility.\end{proof}

\begin{proposition}\label{P1.2a}{\em(Optimizing sequences for $v$ and $v^+$)}
There exists a sequence of stopping times $(\theta^n)_{n \in \mathbb{N}}$
with $\theta^n $ in $ \mathcal{T}_{S}$ (resp.   $\mathcal{T}_{S^+}$)  such that the sequence $(
E[\phi(\theta^n)\, |\,\mathcal{F}_S])_{n \in \mathbb{N}}$ is increasing and such that
\[v(S) \quad {\mbox{\rm{(resp.}  } } v^+(S) \mbox{\rm{)}}\quad = \lim_{n \to \infty} \uparrow E[\phi(\theta^{n}) \, |\,\mathcal{F}_{S}]\quad
\mbox{\rm a.s.}\]
\end{proposition}
\begin{proof}  Again, the arguments are the same for $(v(S), S\in \mathcal{T}_0)$ and  $(v^+(S), S\in \mathcal{T}_0)$.
We prove the property only for  $(v^+(S), S\in \mathcal{T}_0)$.  For each $S$ $\in$ $\mathcal{T}_0$, 
one can show that the set $(E[\phi(\theta)\, |\,\mathcal{F}_S], \; \theta \in \mathcal{T}_{S^+} )$ is closed under pairwise maximization. Indeed, let $\theta,\theta'\in \mathcal{T}_{S^+}$. Put $A=\{\, E[\phi(\theta')\,|\,\mathcal{F}_S]\leq 
 E[\phi(\theta)\,|\,\mathcal{F}_S]\,\}$. One has $A\in\mathcal{F}_S$. Put $\tau=\theta\mathbb{1} _A+\theta' \mathbb{1} _{A^c}$. Then $\tau$ $\in$ $\mathcal{T}_{S^+}$. It is easy to check that $ E[\phi(\tau)\,|\,\mathcal{F}_S]=$ $ E[\phi(\theta)\,|\,\mathcal{F}_S]\vee $ $ E[\phi(\theta')\,|\,\mathcal{F}_S]$.
The result follows by a classical result on essential suprema (Neveu (1975)).
\end{proof} 
An admissible family $( h(\theta) , \; \theta \in \mathcal{T}_{0} )$ is said to be a \emph{supermartingale family} (resp. a \emph{martingale family}) if for any 
$\theta, \theta^{'}$ $ \in$ $\mathcal{T}_0$ such that $\theta \geq \theta^{'}$ a.s., 
\begin{eqnarray*}
E[h(\theta) \, |\,\mathcal{F}_{\theta^{'} }]  \leq  h(\theta^{'}) \quad
\,\mbox{a.s.,} && {\rm (resp.} \quad  E[h(\theta) \, |\,\mathcal{F}_{\theta^{'} }]  =  h(\theta^{'}) \quad
\,\mbox{a.s.).}
\end{eqnarray*}

We now prove that both $v$ and $v^+$ are supermartingale families and that the value function $v$ is characterized as the Snell envelope family associated with the reward $\phi$. More precisely:
\begin{proposition}\label{prop.SuperM} 
The two following properties hold.
\begin{itemize}
\item
The admissible families  $(v(S), S\in \mathcal{T}_0)$ and  $(v^+(S), S\in \mathcal{T}_0)$ are supermartingale families.
\item
The value function family  $(v(S), S\in \mathcal{T}_0)$ is characterized as the \emph{Snell envelope family}  associated with $(\phi(S) , S\in \mathcal{T}_0)$, that is the smallest supermartingale family which is greater (a.s.) than $(\phi(S), S\in \mathcal{T}_0)$
\end{itemize}
\end{proposition}

\begin{proof}  Let us prove the first  point for $v^+$.
Fix $S \geq S^{'}$ a.s.. By Proposition \ref{P1.2a}, there 
exists an optimizing sequence $(\theta^n)$ for $v^+(S)$. By the monotone convergence theorem, 
$\displaystyle{
E[v^+(S) \, |\, \mathcal{F}_{S^{'} }]  =  \lim_{n\to\infty} E[\phi(\theta^{n}) \, |\,\mathcal{F}_{S^{'} }] \,
}$ a.s.. Now, for each $n$, since $\theta^n \in \mathcal{T}_{(S')^+}$, we have  $E[\phi(\theta^n)\,|\,\mathcal{F}_{S'}]\leq v^+(S')$ a.s.
Hence, $E[v^+(S) \, |\, \mathcal{F}_{S^{'} }]   \leq  v^+(S')$ a.s., which gives 
 the supermartingale property of $v^+$. The supermartingale property of $v$ can be proved by using the same arguments.
 
Let us prove the second point (which is classical). First, we clearly have that $(v(S) , S\in \mathcal{T}_0)$ is a supermartingale family and that for each $S \in \mathcal{T}_0$, $v(S) \geq  \phi(S)$ a.s.
Let us prove that is the smallest.  Let $(v'(S),S\in \mathcal{T}_0)$ be a supermartingale family such that for each $\theta \in \mathcal{T}_0$, $v'(\theta)$ $\geq$ $\phi(\theta)$ a.s. Let $S\in \mathcal{T}_0$. By the properties of $v'$, for all $\theta\in \mathcal{T}_S$,
$v'(S)\geq E[v'(\theta)\,|\,\mathcal{F}_S]\geq $  $E[\phi(\theta)\,|\,\mathcal{F}_S]$ a.s. Taking the supremum over $\theta\in \mathcal{T}_S$, we have $v'(S)\geq v(S)$ a.s.
\end{proof} 

The following proposition, known as the optimality criterion, gives a characterization of optimal stopping times for the $v(S)$.

\begin{proposition}\label{prop.criterion}\emph{(Optimality criterion) }
Let $S$ $\in$ $\mathcal{T}_0$ and let $\theta_{*} \in \mathcal{T}_S$ be such that $E[\phi(\theta_{*})]<\infty$.
The following three assertions are equivalent
\begin{enumerate}
\item
$\theta_{*}$ is $S$-optimal for $v(S)$, that is 
\begin{equation}\label{So}
v(S) = E[\phi( \theta_{*}) \, |\,\mathcal{F}_{S}] \quad
\,\mbox{a.s.}
\end{equation}
\item The following equalities hold:  
$v ( \theta_{*}) = \phi( \theta_{*})\quad \,\mbox{a.s.,} \quad  {\rm and} \quad E[v(S)]= E[ v( \theta_{*})].$
\item The following equality holds: $E[v(S)]= E[ \phi( \theta_{*})].$
\end{enumerate}
\end{proposition}
\begin{remark}\label{rem.opt}  
Note that since the value function is a supermartingale family, equality $E[v(S)]= E[ v( \theta_{*})]$ is 
equivalent to the fact that the family $( v(\theta), \theta \in \mathcal{T}_{[S, \theta_{*}]} )$ is a martingale family, that is for all 
$\theta, \theta^{'}$ $\in$ $\mathcal{T}_0$ such that $S \leq \theta, \theta^{'} \leq \theta_{*} $ a.s., $v(\theta)= E[ v( \theta^{'})\,|
\, \mathcal{F}_{\theta}]$ a.s. on $\{\theta \leq \theta^{'}\}$ (which can also be written $\left( v((\theta \vee S) \wedge 
\theta_{*}), \theta \in \mathcal{T}_0\right)$ is a martingale family).
\end{remark}
\begin{proof}  Let us show that 1) implies 2).
Suppose 1) is satisfied. 
Since the value function $v$ is a supermartingale family greater that $\phi$, we have clearly 
\[v(S) \geq E[v(\theta_{*}) \, |\,\mathcal{F}_{S }] \geq E[\phi (\theta_{*}) \, |\,\mathcal{F}_{S }] \;\mbox{a.s. } \]
Since equality (\ref{So}) holds, this implies that the previous inequalities are actually equalities. \\
In particular,  $E[v(\theta_{*}) \, |\,\mathcal{F}_{S }] =  E[\phi (\theta_{*}) \, |\,\mathcal{F}_{S }] $ a.s. but as inequality 
$v(\theta_{*}) \geq \phi(\theta_{*})$ holds a.s., and as $E[\phi(\theta_{*})]<\infty$, we have $v(\theta_{*}) = \phi(\theta_{*})$ a.s..\\
Moreover, $v(S) = E[v(\theta_{*}) \, |\,\mathcal{F}_{S }] $ a.s. which gives $E[v(S)]= E[ v( \theta_{*})]$. Hence, 2) is satisfied. 

Clearly,  2) implies 3). It remains to show that 3) implies 1).
\\
Suppose 3) is satisfied. 
Since $v(S) \geq E[\phi(\theta_{*}) \, |\,\mathcal{F}_{S }] $ a.s., this gives $v(S) = E[\phi(\theta_{*}) \, |\,\mathcal{F}_{S }] $ a.s.. Hence, 1) is safisfied.
\end{proof} 
\begin{remark} \label{eo}
It is clear that 
by 3) of Proposition \ref{prop.criterion}, a stopping time 
$\theta_{*} \in \mathcal{T}_S$ such that $E[\phi(\theta_{*})]<\infty$ is optimal for $v(S)$ if and only if it is optimal for $E[v(S)]$, that is
\[E[v(S)]=\sup_{\theta\in \mathcal{T}_S }E[ \phi( \theta)]=  E[\phi(\theta_{*})].\]
\end{remark}
We state the following property (which corresponds to Proposition D.3 in Karatzas and Shreve (1998)):
\begin{proposition}\label{prop.vv+} 
For all $S\in \mathcal{T}_0$, $v(S)= \phi(S)\vee v^+(S)$ a.s.
\end{proposition}
\begin{proof}  Note first that $v(S)\geq v^+(S)$ a.s. and that $v(S)\geq \phi(S)$ a.s., which yields the inequality $v(S)\geq \phi(S)\vee v^+(S)$ a.s. It remains to show the other inequality. 
Fix $\theta$ $\in$ $\mathcal{T}_S$. First, the following inequality holds:
\begin{equation} \label{eq.premia}
 E[\phi( \theta)\, |\,\mathcal{F}_{S}] \,\mathbb{1}_{\{\theta >S\}} \leq v^+(S) \,\mathbb{1}_{\{\theta >S\}} \,\,\,\mbox{\rm a.s.}
 \end{equation}
 Indeed, since the random variable $\overline \theta$ defined by $\overline \theta := \theta\, \mathbb{1}_{\{\theta >S\}} + T\, \mathbb{1}_{\{\theta \leq S\}}$ belongs to $\mathcal{T}_{S^+}$, one has
 $E[\phi(\overline \theta)\, |\,\mathcal{F}_{S}] \leq v^+(S) $ a.s. and hence 
 \[E[\phi( \theta)\, |\,\mathcal{F}_{S}] \,\mathbb{1}_{\{\theta >S\}} = E[\phi(\overline \theta)\, |\,
\mathcal{F}_{S}]\, \mathbb{1}_{\{\theta >S\}}\leq v^+(S) \,\mathbb{1}_{\{\theta >S\}}\,\,\,\mbox{\rm a.s.}\]
 and thus $
E[\phi(\theta)\,|\,\mathcal{F}_S] = \phi(S)\, \mathbb{1} _{\{\theta =S \}}
+ E[\phi(\theta)\,|\,\mathcal{F}_S]\, \mathbb{1} _{\{\theta >S \}}
\leq \phi(S)\, \mathbb{1} _{\{\theta =S \}}
+ v^+(S)\, \mathbb{1} _{\{\theta >S \}}$ a.s.\,Therefore,
 \[E[\phi(\theta)\,|\,\mathcal{F}_S]\leq \phi(S)\vee v^+(S)\,\,\,{\rm a.s.}\]
 By taking the essential supremum over $\theta \in \mathcal{T}_S$, we derive that $v(S)\leq \phi(S)\vee v^+(S)$ a.s. and the proof is ended.
\end{proof} 
We now provide a useful regularity property for the strict value function family.
\subsection*{Right continuity property of the strict value function}
\begin{definition}\label{def.RCE}
An admissible family $(\phi(\theta),\, \theta\in \mathcal{T}_0)$ is said to be  \emph{right continuous along stopping times in expectation (RCE)} if for any  $\theta\in \mathcal{T}_0$ and for any sequence  of stopping times  $(\theta_n)_{n\in\mathbb{N}}$ such that  $\theta_n\downarrow \theta$  one has
$\displaystyle{E[\phi(\theta)]=\lim_{n \to\infty} E[\phi(\theta_n)]} $.
\end{definition}
The following localization property holds.
\begin{lemma} \label{lemme.phiRCE}
Let $(\phi(\theta),\, \theta\in \mathcal{T}_0)$ be a RCE family. Then, for each $S$ $\in$ $\mathcal{T}_0$ and $A$ $\in$ $\mathcal{F}_S$, the family $(\phi(\theta)\mathbb{1} _A , \; \theta \in \mathcal{T}_S )$ is RCE. 
\end{lemma}
\begin{proof}  Note that if $(\phi(\theta),\, \theta\in \mathcal{T}_0)$ is an admissible family, then for each  $S$ $\in$ $\mathcal{T}_0$ and $A$ $\in$ $\mathcal{F}_S$, the family $(\phi(\theta)\mathbb{1} _A, \theta \in \mathcal{T}_S)$ can easily be shown to be $S$-admissible, that is, to satisfy properties 1) and 2) of Definition \ref{def.admi} with $\mathcal{T}_0$ replaced by $\mathcal{T}_S$.\\
Fix $\theta$ $\in$ $\mathcal{T}_S$.  Let $(\theta_n)_{n \in \mathbb{N}}$ be a nonincreasing sequence of stopping times such that $\theta_n \downarrow \theta$. For each $n$, let $\overline \theta_n := \theta_n \mathbb{1} _A + T \mathbb{1} _{A^c}$ and $\overline \theta := \theta \mathbb{1} _A + T \mathbb{1} _{A^c}$. We clearly have $\overline \theta_n \downarrow \overline \theta$. Hence, since $(\phi(\theta),\, \theta\in \mathcal{T}_0)$ is RCE , it follows that $\lim_{n\to \infty}E[\phi(\overline \theta^n)]= E[\phi( \overline \theta)],$ which clearly yields that $\lim_{n\to \infty}E[\phi(\theta^n)\mathbb{1} _A] = E[\phi( \theta)\mathbb{1} _A].$ 
\end{proof} 
We now show that the strict value function $( v^+(S) , \, S \in \mathcal{T}_0 )$ is RCE (without any regularity assumption on the reward $\phi$). This result is close to Proposition D.3 in Karatzas and Shreve (1998).
\begin{proposition}\label{prop.RCEv+}\emph{(RCE property for $v^+$)}
Let $(\phi(\theta),\, \theta\in \mathcal{T}_0)$ be an admissible family.\\
The associated strict value function family $(v^+(\theta) , \; \theta \in \mathcal{T}_0 )$ is RCE.
\end{proposition}
\begin{remark}\label{rem.RCEA}
Let $S$ $\in$ $\mathcal{T}_0$ and $A$ $\in$ $\mathcal{F}_S$. Since by the previous proposition, $( v^+(\theta) , \theta \in \mathcal{T}_0 )$ is RCE, Lemma \ref{lemme.phiRCE} implies that the family $( v^+(\theta)\mathbb{1} _A , \; \theta \in \mathcal{T}_S )$ is RCE.\\
In particular, the RCE property of $(v^+(\theta)\mathbb{1} _A , \; \theta \in \mathcal{T}_S )$ at $S$ gives that for each non increasing sequence of stopping times $(S_n)_{n \in \mathbb{N}}$  such that $S_n \downarrow S$, we have
\begin{equation*}
\displaystyle{E[v^+(S) \mathbb{1} _A]=\lim_{n \to\infty} E[v^+(S_n)\mathbb{1} _A]}.
\end{equation*}
\end{remark}
\begin{proof} 
Since $(v^+(\theta),\theta \in \mathcal{T}_0)$ is a supermartingale family, the function $\theta\mapsto E[ v^+(\theta)]$ is a non increasing 
function of stopping times. 
 Suppose it is not RCE at $\theta \in \mathcal{T}_0$. We  first consider the case when $E[v^+(\theta)]<\infty$. Then there exists a constant $\alpha>0$ and  a 
 sequence of stopping times $(\theta_n)_{n\in\mathbb{N}}$  such that $\theta_n\downarrow \theta$ and
\begin{equation*}\label{eq.vSna}
\lim_{n\to \infty}\uparrow E[v^+(\theta_n)]+\alpha \leq E[v^+(\theta)].
\end{equation*}
 One can easily show, by using an  optimizing sequence of stopping time for $v^+(\theta)$ (Proposition \ref{P1.2a}) that $\displaystyle{E[v^+(\theta)]= \sup_{\tau \in \mathcal{T}_{S^+}}E[\phi(\tau)]}$.
 Therefore
 there exists $\theta'\in \mathcal{T}_{\theta^+}$ such that 
 \begin{equation}\label{eq.aba}
  \lim_{n\to \infty}\uparrow E[v^+(\theta_n)]+\frac{\alpha}{2} \leq E[\phi(\theta')].
  \end{equation}
Let us first consider the simpler case where $\theta<T$ a.s.\\  In this case,  $\theta'\in \mathcal{T}_{\theta^+}$ implies that $\theta'>\theta$ a.s.; 
one has $\{\theta'>\theta\}= \displaystyle{\bigcup_{n\in \mathbb{N}}} \uparrow \{\theta'>\theta_n\}$ and we have 
$E[\phi(\theta')]=  \displaystyle{\lim_{n\to \infty}}\uparrow E[\mathbb{1}_{\{\theta'>\theta_n\}} \phi(\theta')].$ Hence, there exists $n_0$ such that
\[\lim_{n\to \infty}\uparrow E[v^+(\theta_n)]+\frac{\alpha}{4} \leq E[\mathbb{1}_{\{\theta'>\theta_{n_0}\}}\phi(\theta')].\]
Define the stopping time  $\overline \theta := \theta'\mathbb{1}_{\{\theta'>\theta_{n_0}\}}+T \mathbb{1}_{\{\theta'\leq \theta_{n_0}\}}. $
 One has $\overline \theta > \theta_{n_0}$ a.s. which gives by the positivity of $\phi$ that $ E[\mathbb{1}_{\{\theta'>\theta_{n_0}\}}\phi(\theta')]\leq E[\phi(\overline \theta)]\leq E[v^+(\theta_{n_0})].$ Finally, 
 \begin{equation}\label{eq.prem}
E[v^+(\theta_{n_0})]+\frac{\alpha}{4} \leq \lim_{n\to \infty}\uparrow E[v^+(\theta_n)]+\frac{\alpha}{4} \leq E[v^+(\theta_{n_0})].
\end{equation}
  which gives the expected contradiction.
  
Let us now consider a general  $\theta\in \mathcal{T}_0$.\\ 
Since $\theta'\in \mathcal{T}_{\theta^+}$, we have
$E[\phi(\theta')]= E[\phi(\theta')\mathbb{1}_{\{T>\theta\}}] 
  + E[\phi(T)\mathbb{1}_{\{\theta=T\}}].$
Since, by definition of $\mathcal{T}_{\theta^+}$,  $\theta'>\theta$ a.s. on $\{T>\theta\}$, it follows that  
\[ E[\phi(\theta')\mathbb{1}_{\{T>\theta\}}]=  \displaystyle{\lim_{n\to \infty}}\uparrow E[\mathbb{1}_{\{\theta'>\theta_n\}
\cap \{T>\theta\}} \phi(\theta')].\] This with (\ref{eq.aba}) implies that there exists $n_0$ such that 
\[\lim_{n\to \infty}\uparrow E[v^+(\theta_n)]+\frac{\alpha}{4} \leq E[\mathbb{1}_{\{\theta'>\theta_{n_0}\}\cap \{T>\theta\}}\phi(\theta')]
+ E[\phi(T)\mathbb{1}_{\{\theta=T\}}].\]
 Put $\overline \theta = \theta'\mathbb{1}_{\{\theta'>\theta_{n_0}\}\cap \{T>\theta\}}+T \mathbb{1}_{\{\theta'\leq \theta_{n_0}\}
 \cap \{T>\theta\}} +T \mathbb{1}_{\{T=\theta\}}$. One has $\overline \theta$ $\in$ $\mathcal{T}_{\theta_{n_0}^+}$. Hence,\\
  $ E[\mathbb{1}_{\{\theta'>\theta_{n_0}\}\cap \{T>\theta\}}\phi(\theta')]
 + E[\phi(T)\mathbb{1}_{\{\theta=T\}}] \leq E[\phi(\overline \theta)]\leq E[v^+(\theta_{n_0})]$.
  Finally, we derive again (\ref{eq.prem}) which gives the expected contradiction.
\par
  In the case where $E[v^+(\theta)]= \infty$, by similar arguments, one can show that when  $\theta_n\downarrow \theta$ the limit $\displaystyle{\lim_{n\to \infty}}E[v^+(\theta_n)]$ cannot be finite. The strict value function $( v^+(\theta) , \; \theta \in \mathcal{T}_0)$ is thus RCE.
 \end{proof} 
 We now state a useful lemma.
\begin{lemma}\label{lemme.un}
 Let $(\phi(\theta),\theta\in \mathcal{T}_0)$ be an admissible family.
 For each $\theta, \, S$ $\in$ $\mathcal{T}_0$, we have
\[ E[v(\theta) | \mathcal{F}_{S}] \leq v^+(S)\,\,\,{\rm a.s.} \,\,\, {\rm on} \,\,\,\{\theta > S \}. \]
 \end{lemma}
\begin{proof} 
 Recall that there exists an optimizing sequence of stopping times $(\theta^n)$ with $\theta^n $ in $ \mathcal{T}_{\theta}$  such that 
 $\displaystyle v(\theta)  = \lim_{n \to \infty} \uparrow E[\phi(\theta^{n}) \, |\,\mathcal{F}_{\theta}]\quad
\,\mbox{a.s..}$\\
By taking the conditional expectation, we derive that a.s. on $\{\theta > S \}$,
 \[ E[v(\theta) | \mathcal{F}_{S}] = E[ \lim_{n \to \infty} \uparrow E[\phi(\theta^{n}) \, |\,\mathcal{F}_{\theta}]  |\, \mathcal{F}_{S}]
  =   \lim_{n \to \infty} \uparrow  E[\phi(\theta^{n}) \,  | \,\mathcal{F}_{S}],\]
 where the second equality follows from the monotone convergence theorem for conditional expectation.\\  
 Now, on $\{\theta > S \}$,  since
 $\theta^{n} \geq \theta > S$ a.s., by inequality (\ref{eq.premia}), we have
$ E[\phi (\theta^n) | \mathcal{F}_{S}] \leq v^+(S) $ a.s.
Passing to the limit in $n$ and using the previous equality gives that $E[v(\theta) | \mathcal{F}_{S}] \leq v^+(S)$ a.s. on $\{\theta > S \}$.
\end{proof} 
\begin{proposition}\label{prop.vegalvplus}
Let $(\phi(\theta),\theta\in \mathcal{T}_0)$ be an admissible family of random variables such that $\displaystyle{v(0)= \sup_{\theta\in \mathcal{T}_0}E[ \phi(\theta)]<\infty}$.
Suppose that $( v(S) , \; S \in \mathcal{T}_0 )$ is RCE. Then for each $S \in \mathcal{T}_0$, $v(S)= v^+(S)$ a.s.
\end{proposition}
\begin{proof} 
Fix $S \in \mathcal{T}_0$. For each $n \in \mathbb{N}^*$, put $S_n := (S+ \frac{1}{n}) \wedge T$. Clearly  $S_n \downarrow S$ and for each $n$, $S_n \in \mathcal{T}_{S^+}$ (that is $S_n>S$ a.s. on $\{S<T \}$). 
 By Lemma \ref{lemme.un}, for each $n \in \mathbb{N}$, and a.s. on $\{S<T\}$ we have $ E[v^+(S_n) | \mathcal{F}_{S}] \leq E[v(S_n) | \mathcal{F}_{S}] \leq v^+(S)$. By taking the expectation, we have
 \[E[v^+(S_n) \mathbb{1} _{\{S<T\}} ] \leq E[v(S_n)\mathbb{1} _{\{S<T\}}] \leq E[v^+(S)\mathbb{1} _{\{S<T\}}].\]
  Now, on $\{S=T\}$, for each $n$, $S_n =T$ a.s. and  $v^+(S_n)$ $=$ $v(S_n)$ $=$ $v^+(S)$  $=$ $\phi(T)$ a.s., 
  therefore
 \[E[v^+(S_n)  ] \le E[v(S_n)] \le E[v^+(S)] \]
  which leads, by using the RCE property of $v^+$ to $E[v^+(S)]=E[v(S)]$, but as 
 $v^+(S)\le v(S)$ a.s. and $E[v(S)]\le v(0)<\infty$ we obtain $v(S)=v^+(S)$ a.s.
 \end{proof}

\begin{remark}\label{vRCE} 
Recall that in the particular case where $(\phi(\theta),\theta\in \mathcal{T}_0)$ is supposed to be RCE, the value function $( v(S) , \; S \in \mathcal{T}_0)$ is RCE (see Lemma 2.13 in El Karoui (1981) or Proposition 1.5 in Kobylanski et all (2011)).
\end{remark}

\section{Optimal stopping times}

The main aim of this section is to prove the existence of an optimal stopping time under some minimal assumptions. We stress on that the proof of this result is short and only based on the basic properties shown in the previous sections.

 We use a penalization method as the one introduced by Maingueneau (1978) in the case of a reward process. 
 
More precisely, suppose that $v(0) < \infty$ and fix $S \in \mathcal{T}_0$.
In order to show the existence of an optimal stopping time for $v(S)$, we  first construct for each $\varepsilon$ $\in ]0,1[$, an $\varepsilon$-optimal  stopping time $\theta_{\varepsilon} (S)$ for $v(S)$, that is such that 
\begin{equation*}\label{la}
(1- \varepsilon)v(S)  \, \leq  \,  E [ \phi(\theta_{\varepsilon} (S) )|\mathcal{F}_S ].
\end{equation*}
The existence of an optimal stopping time is then obtained by letting $\varepsilon$ tend to $0$.
 
\subsection{Existence of  epsilon-optimal stopping times}

In the following, in order to simplify notation, we make the change of variable 
$\lambda := 1- \varepsilon$. We now show that if the reward is right upper semicontinuous over stopping times 
in expectation, then, for each $\lambda$ $\in ]0,1[$, there exists an $(1- \lambda)$-optimal stopping time for $v(S)$.
 
Let us now precise the definition of these stopping times.
Let $S$ $\in$ $\mathcal{T}_0$.\\  For $\lambda$ $\in$ $]0,1]$, let us introduce the following $\mathcal{F}_S$-measurable random variable 
\begin{equation}\label{tlpS} \theta^{\lambda} (S):=
{\rm ess} \inf \;\mathbb{T}^{\lambda}_S \quad \mbox{where} \quad \mathbb{T}^{\lambda}_S :=\{\,\theta \in \mathcal{T}_S\, , \,\lambda v(\theta) \leq\phi(\theta) \,\,\mbox {a.s.} \,\}.
\end{equation}

Let us first provide some preliminary properties of these random variables.

\begin{lemma}\label{lemme.tl}\mbox{}\\
 \mbox{}\qquad1.\;  For each $\lambda \in]0,1]$ and each $S\in \mathcal{T}_0$, one has
$\theta^{\lambda} (S)\ge S$ a.s.,\\
 \mbox{}\qquad2.\;  Let $S\in \mathcal{T}_0$ and $\lambda, \lambda' \in]0,1]$. If $\lambda\le \lambda'$,
then $\theta^\lambda(S)\le \theta^{\lambda'}(S)$ a.s.\\
 \mbox{}\qquad 3. For $\lambda \in]0,1]$ and $S,S'\in \mathcal{T}_0$, 
$\theta^{\lambda} (S)\le \theta^{\lambda} (S')$ a.s. on $\{S\le S'\}$. \\
\mbox{}\qquad \; In particular, $\theta^{\lambda} (S)= \theta^{\lambda} (S')$ a.s. on $\{S= S'\}$.
\end{lemma} 
\begin{proof} 
The set ${\mathbb{T}}_S^\lambda$ is clearly stable by pairwise minimization. Therefore, there exists a minimizing sequence $(\theta^n)$ in ${\mathbb{T}}_S^\lambda$ such that  $\theta^n\downarrow\theta^{\lambda} (S)$. In particular, $\theta^{\lambda} (S)$ is a stopping time and $\theta^{\lambda} (S)\geq S$ a.s. 
 
 The second point clearly proceeds from ${\mathbb{T}}^{\lambda'}_S\subset  {\mathbb{T}}^{\lambda}_S$ if $\lambda\le \lambda'$.
 
Let us prove point 3. Let $(\theta_n)_n$ and $(\theta'_n)_n$ be minimizing sequences in ${\mathbb{T}}_S^\lambda$ and ${\mathbb{T}}_{S'}^\lambda$ respectively. Define
$ \tilde \theta_n=\theta'_n\mathbb{1}_{\{S\le S'\}}+ \theta_n\mathbb{1}_{\{S>S'\}}.$
Clearly, $\tilde \theta_n$ is a stopping time in ${\mathbb{T}}_S^\lambda$, hence $\theta^\lambda(S)\le \tilde \theta_n$ a.s., and passing to the limit in $n$ we obtain  $\theta^\lambda(S)\le \theta^\lambda(S')\mathbb{1}_{\{S\le S'\}}+ \theta^\lambda(S)\mathbb{1}_{\{S>S'\}}$ a.s, which gives the expected result. 
 \end{proof} 


Let us now introduce the following definition.

\begin{definition} \label{defr} An admissible family $( \phi(\theta), \theta\in \mathcal{T}_0)$ is said to be \emph{right (resp. left) upper semicontinuous in expectation along stopping times (right (resp. left) USCE)} if for all $\theta\in \mathcal{T}_0$ and for all sequences of stopping times $ (\theta_n)$ such that $\theta^n\downarrow \theta$ (resp. $ \theta^n \uparrow \theta$) 
\begin{equation}\label{usce}
E[\phi(\theta) ]\geq \limsup_{n\to \infty} E[\phi(\theta_n)].
\end{equation} 

An admissible family $( \phi(\theta), \theta\in \mathcal{T}_0)$ is said to be \emph{upper semicontinuous in expectation along stopping times (USCE)} if it is right and left USCE.
\end{definition}

\begin{remark} \label{RUSCE}Note that it is clear that if an admissible family $( \phi(\theta), \theta\in \mathcal{T}_0)$ is right (resp. left) USCE, then, 
for each $S$ $\in$ $\mathcal{T}_0$ and each $A$ $\in $ $\mathcal{F}_S$,  $(\phi(\theta) \mathbb{1} _A, \theta \in \mathcal{T}_S)$ is right (resp. left) USCE. The arguments to show this property are the same as those used in Lemma \ref{lemme.phiRCE}.
\end{remark}

The following Theorem holds:

\begin{theorem}\label{thm.eli}
Suppose the reward $(\phi(\theta), \theta\in \mathcal{T}_0)$ is right USCE and
 $v(0)<\infty$. Let
$S$ in $\mathcal{T}_0$.
For each $\lambda \in ]0,1[$, the stopping time $\theta^{\lambda} (S)$ defined by (\ref{tlpS}) is an $(1- \lambda)$-optimal stopping time for  $v(S)$ that is 
\begin{equation*}
\lambda v(S)  \, \leq  \,  E [ \phi(\theta^{\lambda}(S) )|\mathcal{F}_S ].
\end{equation*}
\end{theorem}

\vspace{0,2cm}

The proof of Theorem \ref{thm.eli}  relies on two lemmas. The first one is the following:

\begin{lemma} \label{stepun}
Suppose  the reward family $( \phi(\theta), \theta\in \mathcal{T}_0)$ is right USCE and
$v(0)<\infty$. 
Then, for each $\lambda \in ]0,1[$, the stopping time $\theta^{\lambda} (S)$ satisfies
\begin{equation*}
\lambda v(\theta^{\lambda} (S) )\leq  \phi(\theta^{\lambda} (S) ) \,\,\,{\rm a.s.}
\end{equation*}
\end{lemma}


\begin{remark}\label{ps} We stress on that the right upper semicontinuity along stopping times \emph{in expectation} of the reward family $\phi$ is sufficient to ensure this key property. The proof relies on the definition of $\theta^{\lambda}(S)$ as an essential infimum of a set of stopping times and on the RCE property of the strict value function family $v^+$.
\end{remark}

\begin{proof}  Let $S \in \mathcal{T}_0$ and $A \in\mathcal{F}_{{\theta^\lambda}(S)}$. In order to simplify notation, let us denote $\theta^\lambda(S)$ by $\theta^\lambda$.

Recall that there exists a minimizing sequence $(\theta^n)$ in ${\mathbb{T}}^\lambda_S$. Hence, $\displaystyle{\theta^\lambda= \lim_{n\to \infty} \downarrow \theta^n}$ and, as $v^+$ 
$\leq$ $v$, we have that for each $n$,
\begin{equation}\label{eq.vn}
\lambda v^+(\theta^n) \leq \lambda v(\theta^n)  \leq \phi(\theta^n)\,\,\,{\rm a.s.}
\end{equation}
 
Note that on $\{v( \theta^\lambda) > \phi( \theta^\lambda) \}$, we have $v( \theta^\lambda)= v^+( \theta^\lambda)$ a.s.\, It follows that
\begin{equation} \label{eq.ran}
\lambda E[ v(\theta^\lambda)\mathbb{1} _A]= \lambda E[ v^+(\theta^\lambda)\mathbb{1}_{\{v( \theta^\lambda) > \phi( \theta^\lambda) \}\cap A}] +\lambda  E[\phi( \theta^\lambda)\mathbb{1}_{\{v( \theta^\lambda) = \phi( \theta^\lambda) \}\cap A}].
\end{equation}
Let us consider the first term of the right member of this inequality and let us now use the RCE property of the strict value function family $v^+$. More precisely, by applying Remark \ref{rem.RCEA} to the stopping time $ \theta^{\lambda}$ and to the set 
$ \{v( \theta^\lambda) > \phi( \theta^\lambda) \}\cap A$, we obtain the following equality

\[ \lambda E[ v^+(\theta^\lambda)\mathbb{1}_{\{v( \theta^\lambda) > \phi( \theta^\lambda) \}\cap A}] = \lambda   \lim_{n
\to \infty}E[v^+(\theta^n)\mathbb{1}_{\{v( \theta^\lambda) >\phi( \theta^\lambda) \}\cap A}].\]
By inequality (\ref{eq.vn}), it follows that 
\[\lambda E[ v^+(\theta^\lambda)\mathbb{1}_{\{v( \theta^\lambda) > \phi( \theta^\lambda) \}\cap A}] 
\leq\limsup_{n\to \infty}E[\phi( \theta^n) \mathbb{1}_{\{v( \theta^\lambda) >\phi( \theta^\lambda) \}\cap A}].\]
Consequently, using equality (\ref{eq.ran}), we derive that 
\begin{eqnarray*} \label{run}
\lambda E[ v(\theta^\lambda)\mathbb{1} _A]&\leq&\limsup_{n\to \infty}E[\phi( \theta^n) \mathbb{1}_{\{v( \theta^\lambda) >\phi( \theta^\lambda) \}\cap A}] + E[\phi( \theta^\lambda)\mathbb{1}_{\{v( \theta^\lambda) = \phi( \theta^\lambda) \}\cap A}]\\
&\leq& \limsup_{n\to \infty}E[\phi(\overline \theta^n)\mathbb{1} _A],
\end{eqnarray*}
where for each $n$, $\overline \theta^n := \theta^n \mathbb{1}_{\{v( \theta^\lambda) > \phi( \theta^\lambda) \}\cap A}  + \theta^\lambda \mathbb{1}_{\{v( \theta^\lambda) = \phi( \theta^\lambda) \}\cap A} + \theta^\lambda \mathbb{1}_{A^c}.$

Note that $(\overline \theta^n)$ is a non increasing sequence of stopping times such that $\overline \theta^n \downarrow \theta^\lambda$. 
Let us now use the right USCE 
assumption on the reward family $\phi$. More precisely, by  Remark \ref{RUSCE}, we have 
\[\lambda E[v(\theta^\lambda)\mathbb{1} _A] \leq \limsup_{n\to \infty}E[\phi(\overline \theta^n)\mathbb{1} _A]\leq  E
[\phi( \theta^\lambda)\mathbb{1} _A].\]
 Hence, the inequality $E\left[\left(\phi( \theta^\lambda)- \lambda v(\theta^\lambda) \right)\mathbb{1} _A\right]\geq 0$ 
 holds for each $A \in\mathcal{F}_{\theta^\lambda}$. By a classical result, it follows that $\phi( \theta^\lambda)- 
 \lambda v(\theta^\lambda) \geq 0$ a.s.\,The proof is thus complete.
\end{proof} 
We now state the second lemma:
\begin{lemma}\label{stepdeux}
Let $(\phi(\theta), \theta\in \mathcal{T}_0)$ be an admissible family with $v(0)<\infty$. For each
 $\lambda \in ]0,1[$ and for each $S\in \mathcal{T}_0$,
\begin{equation}\label{eq.ma}
v(S) =  E [ v(\theta^{\lambda} (S))\, | \, \mathcal{F}_S]\quad \,\mbox{a.s.}
\end{equation}
\end{lemma}
\begin{remark}\label{rem.stepdeux}
Note that equality (\ref{eq.ma}) is equivalent to the martingale property of the family $\displaystyle\left(
 v(\theta), \theta \in \mathcal{T}_{[S, \theta^{\lambda} (S)]} \right)$. In other words,  $(v((\theta \vee S) \wedge \theta^{\lambda} (S)), \theta \in \mathcal{T}_0)$ is a martingale family.
\end{remark}
\begin{proof}  The proof consists to adapt the classical penalization method, introduced by Maingueneau (1978) in the case of a continuous process, to our more general framework. It appears that it is clearer and simpler in the setup of families of random variables than in the setup of processes.
Let us define for each $S$ $\in$ $\mathcal{T}_0$, the random variable
$J_{\lambda}(S) =  E [ v(\theta^{\lambda} (S)) \, |\, \mathcal{F}_S]\,.$
It is sufficient  to show that $J_{\lambda}(S)=v(S)$ a.s.\, Since $( v(S), S \in \mathcal{T}_0 )$ is a supermartingale family and since 
$\theta^{\lambda} (S) \geq S$ a.s., we have that 
\[J_{\lambda}(S) =  E [ v(\theta^{\lambda} (S)) \, |\, \mathcal{F}_S] \leq v(S) \quad \,\mbox{a.s.}\]
It remains to show the reverse inequality. This will be done in two steps.
\\
\emph{Step 1:} Let us show that the family $(  J_{\lambda}(S), S\in \mathcal{T}_0 )$ is a supermartingale family.\\
Fix $ S,S'\in \theta\in \mathcal{T}_0$ such that $S'\geq S$ a.s. We have 
$\theta^{\lambda} (S^{'}) \, \geq\,  \theta^{\lambda} (S)\, \quad \,\mbox{a.s.}$
\par
Hence,
$
E[ J_{\lambda}(S^{'} ) \, |\, \mathcal{F}_S] =    E [ v(\theta^{\lambda} (S^{'})) \, |\, \mathcal{F}_S]
                  $ $=$ $   E\left[  E [ v(\theta^{\lambda} (S^{'})) \, |\, \mathcal{F}_{\theta^{\lambda} (S)} ]   \, |\,\mathcal{F}_S \right]$ a.s.
Now, since $( v(S), S\in \mathcal{T}_0 )$ is a supermartingale family, $E [ v(\theta^{\lambda} (S^{'})) \, |\, \mathcal{F}_{\theta^{\lambda} (S)} ]  \, \leq \, v({\theta^{\lambda} (S)})
$ a.s.
Consequently,
\[E[ J_{\lambda}(S^{'} ) \, |\, \mathcal{F}_S] \leq  E [ v(\theta^{\lambda} (S)) \, |\, \mathcal{F}_S]= J_{\lambda}(S) \quad {\rm a.s.}\]
which ends the proof of step 1.
\\
\emph{Step 2:} 
Let us show that $ \lambda v(S) + (1 - \lambda) J_{\lambda}(S) \, \geq \phi(S)$ a.s. for each $ S \in \mathcal{T}_0$, $\lambda\in]0,1[$.
\\
Fix  $ S \in \mathcal{T}_0$ and $\lambda\in]0,1[$. Let $A:=\{\, \lambda v(S) \leq \phi(S) \,\}$. Let us show that $ \theta^{\lambda} (S) = S $ a.s. on $A$. For this, put $\overline S = S \mathbb{1} _A + T \mathbb{1} _{A^c}$. Note that 
 $\overline S$ $\in$ ${\mathbb{T}}^\lambda_S$. It follows that $ \theta^{\lambda} (S) $ $= {\rm ess} \inf \,{\mathbb{T}}^\lambda_S $ $ \leq$ $\overline S$ a.s. which clearly gives $ \theta^{\lambda} (S) \mathbb{1} _A$
 $ \leq$ $\overline S \,\mathbb{1} _A$ $=$ $S \,\mathbb{1} _A$ a.s. Thus, $\theta^\lambda(S)=S$ a.s.\,on $A$.\\
 Hence, 
$J_{\lambda}(S) \,=\,  E [ v(\theta^{\lambda} (S)) \, |\, \mathcal{F}_S]\,  =\, E [ v(S) \, |\, \mathcal{F}_S] \, = \, v(S) $ a.s.\,on $A$, which yields the inequality
\[ \lambda v(S) + (1 - \lambda) J_{\lambda}(S) = v(S) \geq \phi (S)\quad \,\mbox{a.s.\,\, on\,\,}A.\] 
Furthermore, since $A^c=\{\, \lambda v(S) >\phi( S) \,\}$ and since $J_{\lambda}(S)$ is non negative,
\[ \lambda v(S) + (1 - \lambda) J_{\lambda}(S) \, \geq\,  \lambda v(S)\, \geq \,\phi (S)  \quad \,\mbox{a.s.\,\, on\,\,}A^c.\]
The proof of step 2 is complete.
\par
Note now that, by convex combination, the familly $( \lambda v(S) + (1 - \lambda) J_{\lambda}(S) ,S\in \mathcal{T}_0)$ is a supermartingale family. By step 2, it dominates $(  \phi (S),  S\in \mathcal{T}_0)$. Consequently, by the characterization of $(  v(S), S\in \mathcal{T}_0)$ as the smallest supermartingale 
family which dominates $(  \phi (S),  S\in \mathcal{T}_0)$, we have $ \lambda v(S) + (1 - \lambda) J_{\lambda}(S) \geq v(S) \quad \mbox{\rm a.s.}\,$\\
Hence, $J_{\lambda}(S) \geq v(S)$   a.s. because $v(S) < \infty$ a.s.\, and because $\lambda< 1$ (note that the strict inequality is necessary here). Consequently, for each $S \in \mathcal{T}_0$, $J_{\lambda}(S) = v(S)$   a.s.\,
 The proof of Lemma \ref{stepdeux} is ended. 
 \end{proof} 
\begin{proof}[Proof of Theorem \ref{thm.eli}]
By Lemma \ref{stepdeux}, 
and  Lemma \ref{stepun}
\[\lambda v(S) \,=  \, \lambda E [ v(\theta^{\lambda}(S) )|\mathcal{F}_S] \, \leq  \, E [ \phi(\theta^{\lambda}(S) )|\mathcal{F}_S ].\]
In other words, $\theta^{\lambda}(S)$ is $(1- \lambda)$-optimal for $v(S)$.
\end{proof} 
In the next subsection, under the additional assumption of left USCE property of the reward,  we derive from this theorem that the $(1- \lambda)$-optimal stopping times $\theta^\lambda(S)$ tend to an optimal stopping time  for $v(S)$ as $ \lambda\uparrow 1$.
\subsection{Existence result, minimal optimal stopping times, regularity of the value function family}
\subsubsection{Existence result, minimal optimal stopping times}
\begin{theorem}\label{thm.TAO}\emph{(Existence of an optimal stopping time)}\\
Suppose the reward $(\phi(\theta), \theta\in \mathcal{T}_0)$ is such that  $v(0)<\infty$ and is USCE. Let $S$ $\in$ $\mathcal{T}_0$.\\
The stopping time $ \theta_{*}(S)$ defined by
\begin{equation}\label{E.tao}
 \theta_{*}(S):={\rm ess} \inf\{\theta\in \mathcal{T}_S\,, v(\theta)= \phi(\theta) \; {\rm a.s.} \,\}.
\end{equation}
is the minimal optimal stopping time  for $v(S)$. Moreover, $\theta_{*}(S)= \lim_{\lambda\uparrow 1} \uparrow \theta^\lambda(S)$ a.s.
\end{theorem}
\begin{proof}  
The short proof is based on classical arguments adapted to our framework.
Fix $S$ $\in$ $\mathcal{T}_0$. Since the map $\lambda\mapsto \theta^\lambda(S)$ is non decreasing on $]0,1[$, the random variable $\hat \theta(S)$ defined by 
\begin{equation*}\label{hattheta}
\hat \theta(S):= \lim_{\lambda\uparrow 1} \uparrow \theta^\lambda(S)\,
\end{equation*}
is a stopping time. Let us prove that it is optimal for $v(S)$. By Theorem \ref{thm.eli}, $\lambda E[v(S)] \, \leq  \, E [ \phi(\theta^{\lambda}(S) ) ]$ for each 
$\lambda \in ]0, 1[$. 
Letting $\lambda \uparrow 1$ in this last inequality, and $\phi$ is left USCE, we get $E[v(S  )] \leq E[\phi(\hat{\theta}(S) ) ] $ and hence, $E[v(S)]= E[ \phi(\hat{\theta}(S))]$.
Thanks to the optimality criterion 3) of Proposition \ref{prop.criterion}, $\hat{\theta}(S)$ is optimal for $v(S)$.

Let us now show that $\hat \theta (S) = \theta_{*}(S)$ a.s. and that it is the minimal optimal stopping time.
Note first that $\theta_*(S)=\theta^1(S)$, where $\theta^1(S)$ is the stopping time defined by (\ref{tlpS}) with $\lambda=1$. Now, for each $\lambda\le 1$, $\theta^\lambda(S)\le \theta^1(S)=\theta_*(S)$ a.s.\,Passing to the limit as $\lambda$ tends to $1$, we get $\hat\theta(S)\le \theta_*(S)$. By the optimality criterion, if $\theta\in \mathcal{T}_0$ is optimal for $v(S)$, then $v(\theta) = \phi (\theta)$ a.s. This with the definition of $ \theta_{*}(S)$ leads to $\theta \ge \theta_*(S)$ a.s.\\
It follows that, since $\hat \theta (S) $ is optimal for $v(S)$, we have $\hat \theta (S) \geq \theta_{*}(S)$ a.s.\, Hence, 
$\hat \theta (S)=\theta_{*}(S)$ a.s.\, and it is the minimal optimal stopping time for $v(S)$.
\end{proof} 
\begin{remark}\label{croit}  By Lemma \ref{lemme.tl}and as $\theta_{*}(S)= \theta^{1}(S)$ a.s.\,, we have that for each $S,S^{'}$ $\in$ $\mathcal{T}_0$, $\theta_{*}(S) \leq \theta_{*}(S^{'})$ on 
 $\{S \leq S^{'}\}$. In other words, the map $S \mapsto \theta_{*}(S)$ is non decreasing.
  \end{remark}
\subsubsection{Left continuity property of the value function family}
Note first that, without any assumption on the reward family, the value function is right USCE. 
Indeed, from the supermartingale property of $(v(\theta),\theta\in \mathcal{T}_0)$, we clearly have the following property: for each $S\in \mathcal{T}_0$ and each sequence of stopping times $(S_n)$ such that $S_n \downarrow S$, $\displaystyle{\lim_{n \to \infty} \uparrow E[v(S_n)]\leq E[v(S)]}$.

\vline

Define now the property of \emph{left continuity in expectation along stopping times} (LCE property) similarly to the RCE property (see Definition \ref{def.RCE}) with 
$\theta_n \uparrow \theta$ instead of $\theta_n \downarrow \theta$\,.

Using the monotonicity property of $\theta_{*}$ with respect to stopping times (see Remark \ref{croit}), we derive the following regularity property of the value function:
\begin{proposition}\label{prop.vLCE}
If $(\phi(\theta),\theta\in \mathcal{T}_0)$ is USCE  and $v(0)<\infty$, then $(v(S),S\in \mathcal{T}_0)$ is 
left continuous in expectation along stopping times (LCE). 
\end{proposition}
\begin{proof}  Let $S\in \mathcal{T}_0$ and let $(S_n)$ be a sequence of stopping times such that $
S_n \uparrow S$. Let us show that ${\displaystyle \lim_{n \to \infty} E[v(S_n)]= E[v(S)]}$. First of all, note that for each $n$,  $E[v(S_n)] \geq E[v(S)]$. Hence, 
$\displaystyle{\lim_{n \to \infty} \downarrow E[v(S_n)]\geq E[v(S)]}$.

Suppose now by contradiction that $\displaystyle{\lim_{n \to \infty} \downarrow E[v(S_n)]\neq E[v(S)]}$. Then, there exists $\alpha>0$ such that for all $n$, one has $E[v(S_n)] \geq E[v(S)] +\alpha$.
By Theorem \ref{thm.TAO}, for each $n$, the stopping time $\theta^{*}(S_n)\in \mathcal{T}_{S_n} $ (defined by (\ref{E.tao})) is optimal for 
$v(S_n)$. It follows that for each $n$, $E[\phi(\theta_{*}(S_n))] \geq E[v(S)] +\alpha$.
Now, the sequence of stopping times $(\theta_{*}(S_n))$ is clearly non decreasing.
Let $\overline \theta:= \lim_{n \to \infty} \uparrow \theta_{*}(S_n)$. The random variable 
$\overline \theta$ is clearly a stopping time. Using the USCE property of $\phi$, we obtain 
 \[E[\phi(\overline \theta)] \geq E[v(S)] +\alpha\,.\]
Now, for each $n$, $\theta_{*}(S_n)\geq S_n$ a.s.\, By letting $n$ tend to $\infty$, it clearly follows that $\overline \theta\geq S$ a.s., which provides the expected contradiction. 
\end{proof} 
Consequently, the following corollary holds.
\begin{corollary}\label{coroc}
If $(\phi(\theta),\theta\in \mathcal{T}_0)$ is USCE  and $v(0)<\infty$, then $(v(\theta),\theta\in \mathcal{T}_0)$ is USCE. 
\end{corollary}

\subsection{Maximal optimal stopping times}

\subsubsection{A natural candidate}

Let $( \phi(\theta), \theta\in \mathcal{T}_0)$ be an admissible family and 
$(v(\theta), \theta\in \mathcal{T}_0)$ be the associated value function.\\
Fix $S \in \mathcal{T}_0$, and suppose that  $\theta$ is an optimal stopping time for  $v(S)$, then, as a consequence of the optimality criterion (Remark \ref{rem.opt}), the family  $\displaystyle\left(v(\tau), \tau \in  \mathcal{T}_{[S, \theta]} \right)$ is a martingale family. Consider the set 
\begin{equation*}\label{eq.A}
\mathcal{A}_S= \{\theta \in \mathcal{T}_S, \; \mbox{such that}\, \left(v(\tau), \tau \in  \mathcal{T}_{[S, \theta]} \right) \mbox{is a martingale family}\}
\end{equation*}
A natural candidate for the maximal optimal stopping time for $v(S)$ is  thus the random variable $\check{\theta}(S)$ defined by 
\begin{equation}\label{eq.checkun}
 \check{\theta}(S):={\rm ess} \sup \,\, \mathcal{{A}}_S.
\end{equation}
Note that if $v(0) < \infty$, we clearly have: $\check{\theta}(S) ={\rm ess} \sup\{\theta\in \mathcal{T}_S\,, E[v(\theta)]= E[v(S)]  \,\}$.

\begin{proposition}\label{prop.thetamax}
For each  $S \in \mathcal{T}_0$, the random variable $\check{\theta}(S)$ is a stopping time.
\end{proposition}
This proposition is a clear consequence of  the following lemma
\begin{lemma}\label{lemme.stablemax} 
For each  $S \in \mathcal{T}_0$, the set $\mathcal{{A}}_S$ is stable by pairwise maximization.\\
In particular there exists a nondecreasing sequence $(\theta^n)$ in $\mathcal{A}_S$ such that $\theta^n\uparrow \check{\theta}(S)$. 
\end{lemma}

\begin{proof} 
Let $S \in \mathcal{T}_0$ and $\theta_1$, $\theta_2$ $\in$ $\mathcal{{A}}_S$.
 Let us show that 
$\theta_1 \vee \theta_2$ belongs to $\mathcal{{A}}_S$.
Note that this property is intuitive since if $\displaystyle \left(v(\tau), \, \tau \in \mathcal{T}_{[S, \theta_1]}\right)$ and  $\displaystyle \left(v(\tau), \, \tau \in \mathcal{T}_{[S, \theta_2]}\right)$ are 
martingale families, then it is quite clear that $\displaystyle \left(v(\tau),\,\tau\in \mathcal{T}_{[S, \theta_1\vee \theta_2]}\right)$ is
a martingale family. For the sake of completeness, let us show this property. We have clearly that a.s.
\begin{equation}\label{eqt}
 E[v(\theta_1 \vee \theta_2)\, | \, \mathcal{F}_S]= E[v(\theta_2) \mathbb{1} _{\{\theta_2 > \theta_1\}} \, | \, \mathcal{F}_S] 
 + E[v(\theta_1) \mathbb{1} _{\{\theta_1 \geq  \theta_2\}} \, | \, \mathcal{F}_S].
 \end{equation}
 Since $\theta_2$ $\in$ $\mathcal{{A}}_S$, we have that on $\{\theta_2 > \theta_1\}$, 
 $v(\theta_1) = E[v(\theta_2)| \mathcal{F}_{\theta_1}]$ a.s.\, It follows that\\ 
 $E[v(\theta_2) \mathbb{1} _{\{\theta_2 > \theta_1\}}\, | \, \mathcal{F}_S ]= E[v(\theta_1) \mathbb{1} _{\{\theta_2 > \theta_1\}} \, | \, \mathcal{F}_S]$ a.s.\,  This  with 
 equality (\ref{eqt}) gives that $E[v(\theta_1 \vee \theta_2)\, | \, \mathcal{F}_S]= E[v(\theta_1)\, | \, \mathcal{F}_S]$ a.s.. 
 Now, since $\theta_1$ $\in$ $\mathcal{{A}}_S$, $E[v(\theta_1)\, | \, \mathcal{F}_S]= v(S)$ a.s.. 
 Hence, we have shown that $E[v(\theta_1 \vee \theta_2)\, | \, \mathcal{F}_S]= v(S)$ a.s.\, which gives  that 
 $\theta_1 \vee \theta_2$ $\in$ $\mathcal{{A}}_S$. 
 
The second point of the lemma fellows. In particular, $\check{\theta}(S)$ is a stopping time.
\end{proof} 
\subsubsection{Characterization of the maximal optimal stopping time}
Let $S \in \mathcal{T}_0$. In the sequel, we show that $\check{\theta}(S)$  defined by (\ref{eq.checkun}) is the maximal optimal stopping time for $v(S)$. More precisely,
\begin{theorem}\label{Tsmax}\emph{(Characterization of $\check{\theta}(S)$ as the maximal optimal stopping time)} 
Suppose $(\phi(\theta), \theta\in \mathcal{T}_0)$ is right USCE. Suppose that the associated value function $( v(\theta), \theta\in \mathcal{T}_0)$ is LCE with $v(0)<\infty$.\\
For each $S \in \mathcal{T}_0$, $\check{\theta}(S)$ is the maximal optimal stopping time for $v(S)$.
\end{theorem}
\begin{corollary} If $( \phi(\theta), \,\theta\in \mathcal{T}_0)$ is USCE  and $v(0)<\infty$,  then $\check{\theta}(S)$ is optimal for $v(S)$. 
\end{corollary}
\begin{proof}
By Proposition \ref{prop.vLCE}, the value function $( v(\theta), \theta\in \mathcal{T}_0)$ is LCE, and the Theorem applies.
\end{proof}

\begin{remark}
In the previous works, in the setup of processes,  the maximal optimal stopping time  is given, when the  Snell envelope process $(v_t)$ is a right continuous supermartingale process of class $\mathcal D$,  by using the Doob Meyer decomposition of $(v_t)$ and,  in the general case,  by using the Mertens decomposition of $(v_t)$ (see El Karoui (1981)). Thus  fine results of the General Theory of Processes are needed.

In comparison, our definition of $\check{\theta}(S)$ as an essential supremum of a set of stopping times relies on simpler tools of Probability Theory.
\end{remark}
\begin{proof} [Proof of Theorem \ref{Tsmax}] Fix $S \in \mathcal{T}_0$. To simplify the notation, in the following, the stopping time $\check{\theta}(S)$ will be denoted by $\check{\theta}$. \\
\emph{Step 1}: Let us show that $\check{\theta}$ $\in$ $\mathcal{{A}}_S$.\\
By Lemma \ref{lemme.stablemax}, there exists a nondecreasing sequence  $(\theta^n)$ in $\mathcal{A}_S$ such that $\theta^n\uparrow \check{\theta}$.

\noindent For each $n\in\mathbb{N}$, since $\theta^n$ $\in$ $\mathcal{A}_S$, we have
$E[v( \theta^n)] = E[v(S)].$ Now, as $v$ is LCE, by letting 
$n$ tend to $\infty$ gives $E[v( \check{\theta})] = E[v(S)]$,
and therefore $\check{\theta}$ $\in$ $\mathcal{{A}}_S$.\\
 \emph{Step 2}: Let us now show that $\check{\theta}$ is optimal for $v(S)$.\\
 Let $\lambda \in ]0,1[$. By Lemma \ref{rem.stepdeux}, $\displaystyle\left(v(\tau),\,\tau \in \mathcal{T}_{[\check{\theta}, \theta^{\lambda}(\check{\theta})]}\right)$ is a martingale family. Hence, $\theta^{\lambda}(\check{\theta})$ $\in$ $\mathcal{{A}}_S$. The definition of 
$\check{\theta}$ yields that $\theta^{\lambda}(\check{\theta})= \check{\theta}$ a.s.\,\\
Now, since $\theta^{\lambda}(\check{\theta})$ is $(1- \lambda)$-optimal for $v(\check{\theta})$ and $\phi$ is right USCE, it follows  by Lemma \ref{stepdeux} that 
\[\lambda  v(\check{\theta}) \leq E[ \phi (\theta^{\lambda}(\check{\theta}))| \mathcal{{F}}_{\check \theta}]= \phi(\check{\theta}) \quad \,\mbox{a.s.}\]
Since this inequality holds for each $\lambda \in ]0,1[$, we get $ v(\check{\theta})\leq \phi(\check{\theta}),$  and as $ E[ v(\check{\theta})]\ge E[ \phi(\check{\theta})$],
it follows that $v(\check{\theta})=  \phi(\check{\theta})$ a.s.\,, which implies the optimality of $\check{\theta}$ for $v(S)$.\\
 \emph{Step 3}: Let us show that $\check{\theta}$ is the maximal optimal stopping time for $v(S)$. \\
By Proposition \ref{prop.criterion}, we have that each $\theta$ which is optimal for $v(S)$ belongs to $\mathcal{{A}}_S$ and hence 
is smaller than $\check{\theta}$ (since  $\check{\theta}= {\rm ess} \sup \,\mathcal{{A}}_S)$. This gives step 3. \end{proof} 
\begin{remark} Let $(\phi(\theta),\theta\in \mathcal{T}_0)$ be an admissible family of random variables such that $\displaystyle{v(0)<\infty}$. 
Suppose that $v(0)= v^+(0)$. Then, for each $\theta$ $\in$ $\mathcal{T}_0$, $v(\theta)= v^+(\theta)$ a.s. on $\{\theta < \check \theta (0)\}$. Indeed,   the same arguments as in the proof of Proposition \ref{prop.vegalvplus}    apply to $\displaystyle\left(v(\theta),\,\theta \in \mathcal{T}_{[0,\check \theta (0)[} \right)$, which is RCE (it  is a martingale family).
\par
 By using localization techniques (see below), one can prove more generally that, for each $S, \,\theta \in \mathcal{T}_0 $, 
$v(\theta) = v^+(\theta)\,\,\,{\rm a.s.}\,\,\,{\rm on} \,\,\, \{S \leq \theta < \check \theta (S)\} \cap \{v(S)= v^+(S)\}.$
\end{remark}
\section{Localization and case of equality between the reward and the value function family}
Recall that we have shown that for all $S\in \mathcal{T}_0$, $v(S)=\phi(S)\vee v^+(S)$ a.s. (see Proposition \ref{prop.vv+}). Thus, one can wonder if it possible to have some conditions which ensure that $v(S)= \phi(S)$ almost surely on $\Omega$ (or even locally, that is on a given subset $A$ $\in$ $\mathcal{F}_S$). Thisi s be the object of this section.

We first provide some useful localization properties. 
\subsection{Localization properties}
Let $(\phi(\theta),\theta\in \mathcal{T}_0)$ be an admissible family. Let $S \in \mathcal{T}_0$ and $A$ $\in$ $\mathcal{F}_S$. Let $(v_A(\theta),\theta \in \mathcal{T}_S)$ be the value function associated with the admissible reward $(\phi(\theta)\mathbb{1} _A, \theta \in \mathcal{T}_S)$, defined for each $\theta \in \mathcal{T}_S$ by 
\begin{equation}\label{va}
v_A(\theta)={\rm ess} \sup_{\tau \in \mathcal{T}_{\theta}} E[\phi(\tau)\mathbb{1} _A \, |\,\mathcal{F}_{\theta}]\, ,
\end{equation}
and let $(v_A^+(\theta),\theta \in \mathcal{T}_S)$ be the strict value function associated with the same reward, defined for each $\theta \in \mathcal{T}_S$ by 
\begin{equation}\label{eq.vaplus}
v^+_A(\theta)={\rm ess} \sup_{\tau \in \mathcal{T}_{\theta^+}} E[\phi(\tau)\mathbb{1} _A \, |\,\mathcal{F}_{\theta}]\, .
\end{equation}
Note first that the families $(v_A(\theta),\theta \in \mathcal{T}_S)$ and $(v_A^+(\theta),\theta \in \mathcal{T}_S)$ can easily be shown to be $S$-admissible.

We now state the following localization property:
\begin{proposition} \label{prop.vva}
 Let $\{\phi(\theta),\theta\in \mathcal{T}_0\}$ be an admissible family. Let $\theta\in \mathcal{T}_S$ and let $A$ $\in$ $\mathcal{F}_S$. The value functions $v_A$ and $v_{A^+}$ defined by (\ref{va}) and (\ref{eq.vaplus}) satisfy the following equalities
 \[v_A(\theta) = v(\theta) \mathbb{1} _A\,\,\,\,\,\,{\rm and} \,\,\,\,\,\,v_A^+(\theta) = v^+(\theta) \mathbb{1} _A\,\,\,{\rm a.s.}\]
 \end{proposition}
 
\begin{proof} Thanks to the characterization of the essential supremum (see Neveu (1975)), one can easily show that $v(\theta) \mathbb{1} _A$ coincides a.s. with
 ${\rm ess} \sup_{\tau \in \mathcal{T}_{\theta}} E[\phi(\tau)\mathbb{1} _A \, |\,\mathcal{F}_{\theta}]$, that is $v_A(\theta)$. The proof is the same for the strict value function $v^+$. 
 \end{proof} 
 
 \begin{remark}Let $\theta_{*,A}(S)$ and  $\check \theta_{A}(S)$ be respectively the minimal and the maximal optimal stopping times for $v_A$. One can easily show that $\theta_{*,A}(S)= \theta_{*}(S)$ a.s. on $A$ and $\check \theta_{A}(S)= 
 \check \theta (S)$ a.s. on $A$.
 
 Also, we clearly have that for each $S,S^{'}$ $\in$ $\mathcal{T}_0$, $\theta_{*}(S) \leq \theta_{*}(S^{'})$ on 
 $\{S \leq S^{'}\}$ and $\check \theta(S) \leq \check \theta(S^{'})$ on 
 $\{S \leq S^{'}\}$.
  \end{remark}
 
 \subsection{When does the value function coincide with the reward?}
 
We will now give some local strict martingale conditions on $v$ which ensure the a.s. equality between $v(S)$ and $\phi(S)$ for a given stopping time $S$.
We introduce the following notation:
let $X$, $X'$ be real random variables and let $A \in \mathcal{F}$.\\
We say that $X     \not \equiv  X'$ a.s. on $A$ if $P(\{X \neq X'\} \cap A) \neq 0$.

\begin{definition} Let $u=( u(\theta), \theta\in \mathcal{T}_0)$ be a supermartingale family.
Let $S$ $\in \mathcal{T}_0$ and $A$ $\in$ $\mathcal{F}_{S}$.

The family $u$ is said to be a \emph{martingale family on the right at $S$ on $A$} if there exists $ S^{'} \in \mathcal{T}_0$ with ($S  \leq  S^{'}$ and $S     \not \equiv  S^{'}$) a.s. on $A$ such that $\displaystyle \left( u(\tau), \, \tau \in \mathcal{T}_{[S, {S}^{'}]} \right)$ is a martingale family on $A$.

The family $u$ is said to be a \emph{strict supermartingale family on the right at $S$ on $A$} if it is not a martingale family on the right at $S$ on $A$.
\end{definition}
 
 We now provide a sufficient condition to locally ensure the equality between $v(S)$ and $\phi(S)$ for a given stopping time $S$.
   
    \begin{theorem}\label{thm.vphi}
  Suppose $( \phi (\theta),   \theta \in \mathcal{T}_0)$ is right USCE and such that $v(0) < \infty$. Let $S$ $\in$ $\mathcal{T}_0$ and $A$ $\in$ $\mathcal{F}_{S}$ be such that 
  ($S  \leq  T$ and $S     \not \equiv  T$) a.s. on $A$.\\
 If the value function $( v(\theta),   \theta \in \mathcal{T}_0)$ is a strict supermartingale on the right at $S$ on $A$, then 
 $v(S) = \phi (S)$ a.s.\,on $A$.
  \end{theorem}
 \begin{proof}  Note that, in the case where there exists an optimal stopping time for $v(S)$ and where $A= \Omega$, the above property is clear. Indeed, by assumption, the value function is a strict supermartingale on the right at $S$ on $\Omega$. Also, thanks to the optimality criterion, we derive that $S$ is the only one optimal stopping time for $v(S)$ and hence $v(S) = \phi (S)$ a.s.\\
Let us now consider the general case. 
By Theorem \ref{thm.eli}, for each $\lambda \in ]0,1[$, the stopping time $\theta^{\lambda} (S)$ satisfies:
\begin{equation}\label{labis}
\lambda v(S)  \, \leq  \,  E [ \phi(\theta^{\lambda}(S) )\, |\, \mathcal{F}_{S}  ].
\end{equation}
By Remark \ref{rem.stepdeux}, for each $\lambda \in ]0,1[$, the family $\displaystyle\left( v(\theta), \theta \in \mathcal{T}_{[S, \theta^{\lambda} (S)]} \right)$ is a martingale family. Since $(v(\theta),   \theta \in \mathcal{T}_0)$ is supposed to be a strict supermartingale on the right at $S$ on $A$, it follows that $\theta^{\lambda} (S)= S$ a.s.\,on $A$.
Hence, by inequality (\ref{labis}), we have that for each $\lambda \in ]0,1[$,
\[\lambda v(S) \, \leq  \,  \phi(S ) \quad {\rm a.s.} \; {\rm on} \quad A.\]
By letting $\lambda$ tend to $1$, we derive that $v(S)  \, \leq  \,  \phi(S ) $ a.s. on $A$. Since $v(S) \geq \phi(S)$ a.s.\,, it follows that $v(S) = \phi(S)$ a.s.\,on $A$, which completes the proof.
\end{proof} 
\section{Additional regularity properties of the value function}
We first provide some regularity properties which hold for any supermartingale family.
\subsection{Regularity properties of supermartingale families}\label{resu}
\subsubsection{Left and right limits of supermartingale families along stopping times}\label{lagd}
\begin{definition}\label{def.LL}
Let $S$ $\in$ $\mathcal{T}_0$. An admissible family $(\phi(\theta),\theta\in \mathcal{T}_0)$ is said to be  \emph{left limited along stopping times (LL)} at $S$ if 
there exists an $\mathcal{F}_{S^-}$-measurable random variable $\phi(S^-)$ such that,
for any non decreasing sequence  of stopping times  $(S_n)_{n\in\mathbb{N}}$,
\[\phi(S^-)=\lim_{n \to\infty} \phi(S_n) \,\mbox{a.s. on }A[(S_n)],\]
where $A[(S_n)]= \{S_n \uparrow S \; \mbox{and }S_n<S\,\, \mbox{for all}\,\,n\, \}$.
\end{definition}
Recall some definitions and notation. 
Suppose that $S$ $\in$ $\mathcal{T}_{0^+}$. 

A non decreasing sequence of stopping times $(S_n)_{n \in \mathbb{N}}$ is said to \emph{announce} $S$ on $A$ $\in$ $\mathcal{F}$ if 
\[ S_n \uparrow S \,\mbox{a.s. on }A \; \mbox{and }S_n<S \,\mbox{a.s. on }A .\]

The stopping time $S$ is said to be \emph{accessible on} $A$ if  there exists a non decreasing sequence of stopping times $(S_n)_{n \in \mathbb{N}}$ which announces $S$ on $A$.

The \emph{set of accessibility} of $S$, denoted by $A(S)$ is the union of the sets on which $S$ is accessible. 

Let us recall the following result (Dellacherie and Meyer (1977) Chap IV.80).
\begin{lemma}\label{lemme.D-M}
Let $S$ $\in$ $\mathcal{T}_{0^+}$. There exists a sequence of sets $(A_k)_{k \in \mathbb{N}}$ in $\mathcal{F}_{S^-}$ such that  for each $k$, $S$ is accessible on $A_k$, and $A(S) = \cup_k A_k$ a.s.\,
 \end{lemma}
 It follows that, in Definition \ref{def.LL},  the left limit $\phi(S^-)$ is unique on $A(S)$ and the family $(\phi(S^-)\mathbb{1} _{A(S)}, \,S \in \mathcal{T}_0  )$ is admissible.
\begin{theorem} \label{thm.sull}
A supermartingale family $(u(\theta), \theta\in \mathcal{T}_0)$ is left limited along stopping times (LL) at each time $S$ $\in$ $\mathcal{T}_{0^+}$.\\
If $u(0) < +\infty$, then $u(S^-)\mathbb{1} _{A(S)}$ is integrable.
\end{theorem} 
This result clearly follows from the result of Dellacherie and Meyer (1977) quoted above together with the following lemma.
\begin{lemma} \label{sulli}
Let $(u(\theta), \theta\in \mathcal{T}_0)$ be a supermartingale family. Let $S$ be a stopping time in $\mathcal{T}_{0^+}$. Suppose that $S$ is accessible on a measurable subset $A$ of $\Omega$.\\
There exists an $\mathcal{F}_{S^-}$-measurable random variable $u(S^-)$, unique on $A$ (up to the equality a.s.\,), such that, for any non decreasing sequence $(S_n)_{n\in \mathbb{N}}$ announcing $S$ on $A$, one has
\[u(S^-)=\lim_{n\to\infty} u(S_n) \quad \,\mbox{a.s. on } A.\]
If $u(0) < +\infty$, then $u(S^-)\mathbb{1} _A$ is integrable.
\end{lemma} 
\begin{proof}   Let $S$ be stopping time accessible on a set $A$ $\in$ $\mathcal{F}$ and let $(S_n)$ be a sequence announcing $S$ on $A$. 
It is clear that $(u( S_n))_{n\in \mathbb{N}}$ is a discrete non negative supermartingale relatively to the filtration $(\mathcal{F}_{S_n})_{n\in \mathbb{N}}$. By the well-known convergence theorem for discrete supermartingales, there exists a random variable $Z$ such that  
$(u( S_n))_{n\in \mathbb{N}}$ converges a.s. to  $Z$. If $u(0) < +\infty$, then $Z$ is integrable.
Set $u(S^-) := Z$. 

It remains to show that this limit, on $A$, does not depend on the sequence $(S_n)$. 
Let $(S'_n)$ be a sequence announcing $S$ on $A$. Again, by the supermartingales convergence theorem, there exists a  random variable $Z^{'}$ such that  
$(u(S'_n))_{n\in \mathbb{N}}$ converges a.s. to  $Z^{'}$. We will now prove that $Z= Z^{'}$ a.s.\,on $A$. 

For each $n$ and each $\omega$, consider the reordered terms $S^{(0)}(\omega)$ $\le$ $S^{(1)}(\omega)$ 
$\le$ $\cdots$ $\le$ $S^{(2n)}(\omega)$ of the sequence $S_0(\omega)$, $\cdots$, $S_n(\omega)$, $S'_0(\omega)$, $\cdots$, 
$S'_n(\omega)$ and define $\tilde S_n(\omega):= S^{(n)}(\omega)$. It is easy to see that for each $n$, $\tilde S_n$ is a stopping time and that the sequence $(\tilde S_n)$ announces $S$ on $A$. Again, by the supermartingales convergence theorem, there exists a  random variable $\tilde Z$ such that  
$(u(\tilde S_n))_{n\in \mathbb{N}}$ converges a.s. to  $\tilde Z$. Let us show that 
\[Z= \tilde Z\,\,\, \,\,\,{\rm a.s.}\,\,\,{\rm on}\,\,\,A.\]
For almost every $\omega$ $\in$ $A$, as $S_n(\omega)<S(\omega)$ and $S'_n(\omega)<S(\omega)$ for all $n$, the sequence $(\tilde S_n(\omega))$ describes all the values taken by both the sequences $(S_n(\omega))$ and 
$(S'_n(\omega))$ on $A$. Hence, by construction, for each $k$, 
\begin{equation}\label{nk}
A = \cup_{n \geq k} \{\tilde S_n = S_k \} \cap A  
 \end{equation}
almost surely. Without loss of generality, we can suppose that this equality is satisfied everywhere.
Also, by the admissible property of the value function, for each $k$, $n$ $\in$ $\mathbb{N}$, the following equality 
\begin{equation}\label{kn}
 u(S_k)= u( \tilde S_n) \,\,\, \,\,{\rm on}\,\,\, \{\tilde S_n = S_k \}
 \end{equation}
 holds almost surely.
Again, without loss of generality, we can suppose that for each $k$, $n$ $\in$ $\mathbb{N}$, this equality is satisfied everywhere on the set $\{\tilde S_n = S_k \}$. Also, we can suppose that the sequences 
$(u(S_n))$ and $(u(\tilde S_n))$ converge to $Z$ and $\tilde Z$ everywhere on $\Omega$. 

Let $\varepsilon >0$ and $\omega$ $\in$ $A$.\\
Suppose that $Z(\omega)$ and $\tilde Z(\omega)$ are finite. There exists $k_0(\omega)$ $\in$ $\mathbb{N}$ such that for each $n, k$ $\geq$ $k_0(\omega)$, 
\begin{eqnarray}\label{tildez}
\vert u(S_k)(\omega) - Z(\omega) \vert \leq \varepsilon &{\rm and}& \vert u(\tilde S_n)(\omega) - \tilde Z(\omega) \vert \leq \varepsilon. 
\end{eqnarray}
Now, by (\ref{nk}), there exists $n_0(\omega) \geq k_0(\omega)$ such that 
$S_{k_0(\omega)} (\omega) = \tilde S_{n_0(\omega)} (\omega)$. Hence, by (\ref{kn}), 
\[u(S_{k_0(\omega)} )(\omega) = u( \tilde S_{n_0(\omega)} )(\omega).\]
By inequalities (\ref{tildez}), it follows that $\vert  Z(\omega) - \tilde Z(\omega) \vert \leq 2\varepsilon$. Since this inequality holds for each $\varepsilon >0$, we have $Z(\omega) = \tilde Z(\omega)$. Similar arguments show that if $Z(\omega)$ or $\tilde Z(\omega)$ is not finite, then both are not finite. We thus have proven that $Z= \tilde Z$ a.s.\,on $A$. By symmetry, 
$Z^{'}= \tilde Z$ a.s.\,on $A$, which yields the equality $Z^{'}= Z$ a.s.\,on $A$.\\
We have thus shown that $u(S^-)(= Z)$, on $A$, does not depend on the sequence $(S_n)$.   

It remains to show that $u(S^-)$ can be chosen $\mathcal{F}_{S^-}$-measurable. Indeed, the above part of the proof still holds with $A$ replaced by $A[(S_n)]= \{S_n \uparrow S \; \mbox{and }S_n<S\,\, \mbox{for all}\,n\, \}$, which contains $A$, and with $u(S^-)(= Z)$ replaced by $Z\,\mathbb{1} _{A[(S_n)]}$. Note now that $A[(S_n)]$ $\in$ $\mathcal{F}_{S^-}$. Indeed, 
$A[(S_n)] = (\cap_n\{S_n < S\}) \backslash \{\lim S_n <S\}$.
The proof of the lemma is thus complete.
\end{proof} 
\begin{definition}
Let $S$ $\in$ $\mathcal{T}_0$. An admissible family $(\phi(\theta),\theta\in \mathcal{T}_0)$ is said to be  \emph{right limited along stopping times} (RL) at $S$ 
if 
there exists an $\mathcal{F}_{S}$-measurable random variable $\phi(S^+)$ such that,
for any non increasing sequence  of stopping times  $(S_n)_{n\in\mathbb{N}}$,
such that  $S_n\downarrow S$ and $S_n >S$ for each $n$, one has
$\displaystyle{\phi(S^+)=\lim_{n \to\infty} \phi(S_n)} $.
\end{definition}
\begin{theorem} \label{surl}
A  supermartingale family $(u(\theta), \theta\in \mathcal{T}_0)$, with $u(0) < + \infty$, is right limited along stopping times (RL)
at any stopping time $S$ $\in \mathcal{T}_0$. 
\end{theorem} 
\begin{proof} Let $(S_n)_{n\in \mathbb{N}}$ in $T^+_S$ such that $S_n \downarrow S$.\\
Set $Z_n := u(S_{-n})$ and $\mathcal{G}_n:= \mathcal{F}_{S_{-n}}$ for each $n \leq 0$. The sequence $(Z_n)_{n \leq 0}$ is a supermartingale with respect to the non decreasing filtration $(\mathcal{G}_n)_{n\leq 0}$.\\
By a convergence theorem for discrete supermartingales indexed by non positive integers, and uniformly bounded in $L^1$ (see chap.V, Thm.30 in Dellacherie and Meyer (1980)), there exists an integrable random variable $Z$ such that 
the sequence $(Z_n)_{n \leq 0}$ converges a.s. and in $L^1$ to $Z$. We then define $u(S^+)$ by $u(S^+):= Z$.\\ 
It remains to show that this limit does not depend on the sequence $(S_n)$. The proof is not detailed since it is similar to that of the previous theorem. 
\end{proof} 
\subsubsection{Jumps of supermartingale families}
\begin{definition}\label{def.RC-LC}
An admissible family $(\phi(\theta),\theta\in \mathcal{T}_0)$ is said to be \emph{right continuous along stopping times (RC)} if for any  $\theta\in \mathcal{T}_0$ and for any sequence  of stopping times  $(\theta_n)_{n\in\mathbb{N}}$ such that  $\theta_n\downarrow \theta$  one has $ \displaystyle{\phi(\theta)=\lim_{n \to\infty} \phi(\theta_n)} $ a.s.
\end{definition}
The \emph{left continuity along stopping times (LC)} property is defined in a similar way. 
\begin{proposition}\label{droite}
Let $(u(\theta), \theta\in \mathcal{T}_0)$ be a uniformly integrable supermartingale family . 
\begin{itemize}
\item
For each $S$ $\in $ $\mathcal{T}_{0}$, we have $u(S^+)\leq  u(S)$ a.s.
\item Suppose that $(u(\theta),\theta\in \mathcal{T}_0)$ is RCE.
Then, $(u(\theta),\theta\in \mathcal{T}_0)$ is RC, that is, for each $S$ in $\mathcal{T}_0$, $u(S^+) = u(S)$ a.s.
\item
If $(u(\theta),\theta\in \mathcal{T}_0)$ is a martingale family, then $(u(\theta),\theta\in \mathcal{T}_0)$ is RC.
\end{itemize}
\end{proposition}
\begin{proof} 
Let $(S_n)_{n\in \mathbb{N}}$ in $\mathcal{T}_{S^+}$ such that $S_n \downarrow S$.
\par
Let us prove the first point. Thanks to the RL property of $u$, we have
$u(S^+) = \lim_{n\to\infty} u(S_n) $ a.s.\,The supermartingale property of $u$ yields that $E[u(S_n)\,|\,\mathcal{F}_{S}] \leq u(S)$ a.s. for each $n$. By letting $n$ tend to $\infty$ and using the uniform integrability property of $(u(S_n))$, we get $E[u(S^+)\,|\,\mathcal{F}_{S}] \leq u(S)$ a.s.\,Since $u(S^+)$ is $\mathcal{F}_{S}$-measurable, we get $u(S^+)\leq $ $u(S)$ a.s.\,
\par
Let us now prove the second point. Thanks to the RL property of $u$ and the uniform integrability property of $(u(S_n))$, we have
$E[u(S^+)] = \lim_{n\to\infty} E[u(S_n)] = E[u(S)]$, where the last equality follows from the RCE property of $u$.
Now, by the first point, we have $u(S^+)\leq $ $u(S)$ a.s. This with the previous equality leads to the desired result.
\par
The last point is clear.
\end{proof} 
\begin{proposition}\label{prop.gauche}
Let $(u(\theta),\theta\in \mathcal{T}_0)$ be a supermartingale family with $u(0) < + \infty$.\\
Let $S$ $\in $ $\mathcal{T}_{0^+}$ and $(S_n)$ be a non decreasing sequence in $\mathcal{T}_0$ such that $S_n \uparrow S$ a.s.
\begin{itemize} 
\item
We have
\[u(S^-)\geq E[u(S)\,|\,\mathcal{F}_{S^-}] \, \,\mbox{a.s.  on  } \{S_n<S,\,\, \mbox{for all}\,\,n\, \}.\]
\item If $(u(\theta),\theta\in \mathcal{T}_0)$ is LCE, then
\[u(S^-)=E[u(S)\,|\,\mathcal{F}_{S^-}] \quad \,\mbox{a.s.\,\,\,on  }\,\, \{S_n<S,\,\, \mbox{for all}\,\,n\, \}.\]
\end{itemize}
Let $S$ be a predictable stopping time. Then, $u(S^-)\geq E[u(S)\,|\,\mathcal{F}_{S^-}]$ a.s.\,, and if $u$ is LCE, this inequality is an equality.
\end{proposition}
\begin{proof} Let $S$ $\in $ $\mathcal{T}_{0^+}$ and let $(S_n)$ be a non decreasing sequence in $\mathcal{T}_0$ such that $S_n \leq S$ a.s.\,\\
Let us prove the first assertion.Thanks to the LL property of the supermartingale family, we have $\lim_{n\to\infty} u(S_n)= u(S^-)$ a.s.\,on $A[(S_n)_{n\in \mathbb{N}}]=\{S_n \uparrow S \; \mbox{and }S_n<S,\,\, \mbox{for all}\,\,n\, \}$.
Now, since $u$ is a supermartingale family, for each $n$, $u(S_n)\geq E[u(S)\,|\,\mathcal{F}_{S_n}]$ a.s. By letting $n$ tend to $\infty$, we get $\lim_{n \to \infty} u(S_n) \geq E[u(S)\,|\,  \vee_n\mathcal{F}_{S_n}]$ a.s.\, which provides that 
\[u(S^-)\geq E[u(S)\,|\, \vee_n\mathcal{F}_{S_n}] \,\,\, \,\mbox{a.s.  on  }\,\,\, A[(S_n)].\]

From now on, we also suppose that $\lim_{n \to \infty} S_n = S$ a.s.\,\\
We thus have $A[(S_n)]$ $=$ $\{S_n<S,\,\, \mbox{for all}\,\,n\, \}$ and hence belongs to 
$\vee_n\mathcal{F}_{S_n} \cap\mathcal{F}_{S^-}$. By using Lemma \ref{measurability}, it follows that
$E[u(S)\,|\, \vee_n\mathcal{F}_{S_n}]= E[u(S)\,|\,\mathcal{F}_{S^-}]$ a.s. on $ A[(S_n)]$.

It remains to show the second assertion.
Since for almost every $\omega$ $\in$ $ A[(S_n)]^c$, the sequence $(S_n(\omega))$ is constant from a certain rank and thanks to the admissibility property of $u$, it follows that the sequence $(u(S_n)(\omega))$ is also constant from a certain rank. Hence, $\lim_{n\to\infty} u(S_n)= u(S)$ a.s.\,on $ A[(S_n)]^c$. By Fatou's lemma, we have
$E[u(S^-)\mathbb{1} _{A[(S_n)]}] + E[u(S)\mathbb{1} _{A[(S_n)]^c}] \leq \lim_{n\to\infty} E[u(S_n)]$.\\
The LCE property of $u$ yields that $\lim_{n\to\infty} E[u(S_n)]= E[u(S)]$. Hence, 
\begin{equation}\label{eq.inegau}
 E[u(S^-)\mathbb{1} _{A[(S_n)]}] \leq E[u(S)\mathbb{1} _{A[(S_n)]}].
 \end{equation}
Now, since $u$ is a supermartingale family, $u(S^-)\geq E[u(S)\,|\,\mathcal{F}_{S^-}]$ a.s.\,on $A[(S_n)]$. This with inequality (\ref{eq.inegau}) leads to  
$u(S^-)= E[u(S)\,|\,\mathcal{F}_{S^-}]$ a.s.\,on $A[(S_n)]$, which yields the second assertion.
\end{proof} 

Note that, for each admissible family $\left(\phi(\theta), \theta \in \mathcal{T}_0 \right)$, the family of random variables $(E[\phi(S)\,|\,\mathcal{F}_{S^-}], S \in \mathcal{T}_0^p)$, where $\mathcal{T}_0^p$ is the set of predictable stopping times smaller than $T$, can be shown to be admissible. 

\begin{remark} \label{Meyeun}
Let us consider the case when the family $\phi$ is defined via a progressive process $(\phi_t)$. For each $S \in \mathcal{T}_0^p$, we have $E[\phi_S\,|\,\mathcal{F}_{S^-}]$ $=$ $^p\phi_S$ a.s.\, where $(^p\phi_t)$ is the predictable projection of the process $(\phi_t)$ (see 
Th.43 Chap VI  in Dellacherie and Meyer (1980)). Note that the notion of accessible projection of a given progressive process is much more involved. 
\end{remark}
Let $(u(\theta),\theta\in \mathcal{T}_0)$ be a supermartingale family with $u(0) < + \infty$.
For each $S$ $\in $ $\mathcal{T}_{0}$, the \emph{left jump} of $u$ at $S$ is defined by $\Delta u(S) = u(S)- u(S^-)$ on $A(S)$, the set of accessibility of $S$.
We have
\begin{equation}\label{eq.mey}
\Delta u(S)= \left(u(S) - E[u(S)\,|\,\mathcal{F}_{S^-}] \right)+ \left(E[u(S)\,|\,\mathcal{F}_{S^-}]- u(S^-)\right).
\end{equation}
Suppose that $S$ is predictable.\\
Then, by the above proposition, the second term of the left hand side part of equality (\ref{eq.mey}) is non positive,
 and in the particular case when $u$ is LCE, it is equal to $0$.

\begin{remark} \label{Meyedeux} Let us consider the case when the family $u$ is defined via a supermartingale process $(u_t)$. The last assertion gives that, if $S$ is predictable, 
\[u_{S^-}- E[u_{S}\,|\,\mathcal{F}_{S^-}]= u_{S^-}- \,^pu_S \geq 0 \quad  {\rm a.s.}\,,\]
which corresponds to Th. 14 Chap VI in 
Dellacherie and Meyer (1980). This inequality is linked to the jumps of the predictable non decreasing process $(A_t)$ associated to the decomposition of $(u_t)$ (see equality \ref{Meyetrois}).
\end{remark}

Consider now a general stopping time $S$ $\in$ $\mathcal{T}_{0^+}$. Recall that there exists a sequence of sets $(A_k)_{k \in \mathbb{N}}$ in $\mathcal{F}_{S^-}$ such that  for each $k$, $S$ is accessible on $A_k$, and $A(S) = \cup_k A_k$ a.s.\,
 One can easily show that for each $k$, there exists a predictable stopping time $\tau_k$ such that $S= \tau_k$ on $A_k$ a.s.\,(see for example  Lemma 4.7 in \cite{KQ}). It follows that for each $S$ $\in $ $\mathcal{T}_{0^+}$, $\Delta u(S)= \Delta u(\tau_k)$ on $A_k$ a.s.
 
From this, we derive the following property.
\begin{proposition}\label{gauchebis}
Suppose the filtration is left quasicontinuous. 
Let $(u(\theta),\theta\in \mathcal{T}_0)$ be a supermartingale family with $u(0)<+ \infty$. 
\begin{itemize}
\item Suppose that $(u(\theta),\theta\in \mathcal{T}_0)$ is LCE.
Then, $(u(\theta),\theta\in \mathcal{T}_0)$ is LC, that is, for each $S$ in $\mathcal{T}_{0^+}$, $\Delta u(S) =0$ on $A(S)$ a.s.
\item
If $(u(\theta),\theta\in \mathcal{T}_0)$ is a martingale family, then $(u(\theta),\theta\in \mathcal{T}_0)$ is LC.
\end{itemize}
\end{proposition}
\begin{remark} 
When $u$ is defined via a martingale process $(u_t)$, the last assertion implies that, if the filtration is left quasicontinuous, the martingale $(u_t)$ has only totally inaccessible jumps.
\end{remark}

All the above properties hold for the value functions families $v$ and $v^+$ since they are supermartingale families.

\subsection{Complementary properties of the value function}

First, the value functions families $v$ and $v^+$ satisfy the following property.
\begin{proposition}\label{vplusplus}
Let $(\phi(\theta),\theta\in \mathcal{T}_0)$ be a uniformly integrable admissible family. 
For each $S$ $\in$ $\mathcal{T}_0$,
\[v^+(S) = v(S^+)\,\,\,{\rm a.s.}\]
\end{proposition}
\begin{proof}  First, by the second point of Proposition \ref{droite}, and since $v^+$ is RCE, $v^+$ is RC.\\
Let $(S_n)_{n\in \mathbb{N}}$ in $\mathcal{T}_{S^+}$ such that $S_n \downarrow S$. One has $v(S_n) \geq v^+(S_n)$ a.s. for each $n$. Passing to the limit, we have $v(S^+) \geq v^+(S)$ a.s.\,\\
Also, since $S_n > S$ a.s.\,, by Lemma \ref{lemme.un}, we have 
$v^+(S)\geq E[v(S_n)\,|\,\mathcal{F}_{S}]$ a.s. for each $n$. Letting $n$ tend to $+ \infty$, the uniform integrability property of $(v(S_n))$ yields that 
$v^+(S)\geq E[v(S^+)\,|\,\mathcal{F}_{S}]$ $=$ $v(S^+)$ a.s.\,Hence, $v^+(S) = v(S^+)$ a.s.
\end{proof} 
We now provide some local properties of the value function at a stopping time on the left.

First, in the case where the reward is supposed to be USCE, we have the following property.
\begin{proposition} \label{contgauche} Let $(\phi(\theta),\theta\in \mathcal{T}_0)$ be an USCE family such that $v(0)<\infty$.
Then, the associated value function family $(v(\theta),\theta\in \mathcal{T}_0)$ satisfies:\\
For each $S$ $\in $ $\mathcal{T}_{0^+}$ and for each non decreasing sequence $(S_n)$ in $\mathcal{T}_0$ such that $S_n \uparrow S$ a.s.\,, we have
\[v(S^-)=E[v(S)\,|\,\mathcal{F}_{S^-}] \quad \,\mbox{a.s.\,\,\,on  }\,\, \{S_n<S,\,\, \mbox{for all}\,\,n\, \}.\]
Moreover, for each predictable stopping time $S$ in $\mathcal{T}_{0^+}$,
\[v(S^-)=E[v(S)\,|\,\mathcal{F}_{S^-}] \quad \,\mbox{a.s.  }.\]
\end{proposition}

\begin{proof} 
 By Property \ref{prop.vLCE}, the value function $v$ is LCE. Thanks to Proposition \ref{prop.gauche} applied to the supermartingale family $(v(\theta),\theta\in \mathcal{T}_0)$, the result follows.
\end{proof} 

Moreover, in the case where the reward is only supposed to be right USCE, we provide the following result.
 
\begin{theorem} \label{agauche} Suppose $(\phi (\theta),   \theta \in \mathcal{T}_0)$ is a right USCE admissible family and that $v(0) < + \infty$. 
Let $\theta$ $\in$ $\mathcal{T}_0$. Let $(\theta_n)_{n \in \mathbb{N}}$ in $\mathcal{T}_0$ such that $\theta_n \uparrow \theta$.

Suppose that the event $A= A[(\theta_n)]: = \{\theta_n<\theta,\,\, \mbox{for all}\,\,n\, \}$ is non empty. Then, we have
\begin{equation*}
 v(\theta^-)= \overline \phi_C (\theta)\,\,\, \,\,\,  {\rm a.s. }\,\,\,{\rm on}\,\,\,\,\,\, C:= \{E [ v({\theta})  \, |\, \mathcal{F}_{\theta^-} ] < v({\theta^-})\} \cap A,
 \end{equation*}
 where
 \begin{equation}\label{env}
 \displaystyle{
\overline \phi_C (\theta)  :={\rm ess}\!\!\! \!\!\!  \!  \sup_{\!\!\!  \;(\theta_n^{'}) \in \mathcal{A}( \theta, C )} \limsup _{n \to \infty} \phi (\theta_n^{'})\mathbb{1} _{C}},
\end{equation}
and $ \mathcal{A}( \theta, C )$ is the set of non decreasing sequences in $\mathcal{T}_0$ which announce $\theta$ on $C$.
 
Moreover if $( \phi (\theta),   \theta \in \mathcal{T}_0)$ is also LL at $\theta$ on $C$, we have
$ v(\theta^-)= \phi(\theta^-)$ a.s. on $C$.

 \end{theorem}

If $\theta$ is predictable, then for each sequence $(\theta_n)_{n \in \mathbb{N}}$ in $\mathcal{T}_0$ announcing $\theta$, $A[(\theta_n)]= \Omega$ a.s.\,
Hence, $C=B$, where $B:= \{E [ v({\theta})  \, |\, \mathcal{F}_{\theta^-} ] < v({\theta^-})\}$ and thus does not depend on the given sequence $(\theta_n)$ which announce $\theta$, and similarly for $\overline \phi_B (\theta)$. The more tractable property thus follows.

\begin{corollary}\label{predictable}
Suppose $(\phi (\theta),   \theta \in \mathcal{T}_0)$ is a right USCE admissible family and that $v(0)<  \infty$. 
Let $\theta$ in $\mathcal{T}_0$ be a predictable stopping time. We have
\begin{equation*}
 v(\theta^-)= \overline \phi_B (\theta)\,\,\, \,\,\,  {\rm a.s. }\,\,\,{\rm on}\,\,\,\,\,\, B:= \{E [ v({\theta})  \, |\, \mathcal{F}_{\theta^-} ] < v({\theta^-})\} ,
 \end{equation*}
 where $\overline \phi_B (\theta)$ is defined by (\ref{env}).
 
Moreover if $( \phi (\theta),   \theta \in \mathcal{T}_0)$ is also LL at $\theta$ on $B$, we have
$ v(\theta^-)= \phi(\theta^-)$ a.s. on $B$.

\end{corollary}

We stress the importance of this corollary, which allows us to compute the jumps $\Delta A^d_t$ of the predictable non decreasing process $(A_t)$ associated to the decomposition of $(u_t)$ (see Proposition \ref{prop.flatdeuxbis}).
\begin{proof}[Proof of Theorem \ref{agauche}]
Let $B:=\{E [ v({\theta})  \, |\, \mathcal{F}_{\theta^-} ] < v({\theta^-})\}$. Note that $C= B \cap A$.\\
Also, $A$ and $B$ belong to $\mathcal{F}_{\theta^-}$. Hence, $\overline{\phi}_C(\theta)$ is clearly well defined and $\mathcal{F}_{\theta^-}$-measurabme.
\par

We also have that $A\in \mathcal{F}_{\theta^-}\cap (\vee_n \mathcal{F}_{\theta_n}).$

Let us first show that for each $p$ and for each $\lambda\in [0,1[$, $\theta^\lambda(\theta_p)<\theta$ a.s. on $B\cap A$. For this, it is sufficient to
 show that for each $p$, $B \cap A \cap \{\theta^\lambda(\theta_p)\ge \theta\} = \emptyset$ a.s. \\
 Note first that $\{\theta^\lambda(\theta_p)\ge \theta\}=\cap_q\theta^\lambda(\theta_p)\ge \theta_q\}$.
 Hence, $\{\theta^\lambda(\theta_p)\ge \theta\}\in \mathcal{F}_{\theta^-} \cap \vee_n\mathcal{F}_{\theta_n}$. Also, for each 
 $q \ge p$, $E[v(\theta\lambda(\theta_p))| \mathcal{F}_{\theta_q}] = v(\theta_q)$ a.s. on ${\theta\lambda(\theta_p) \ge  \theta_q}$ and hence on ${\theta\lambda(\theta_p) \ge  \theta}$, because $ \{v(\tau ),\tau  \in  \mathcal{T}_{[\theta_p,\theta^\lambda(\theta_p)]}\}$ is a martingale.
Hence, by letting $q$ tend to $\infty$,
\begin{equation}\label{eq-vtmoins}
E[v(\theta^\lambda(\theta_p)| \vee_n \mathcal{F}_{\theta_n}]=v(\theta^-) \quad \mbox{ a.s. on } \quad \{\theta^\lambda(\theta_p)\ge \theta\} \cap A.
\end{equation}   
Now, by a measurability property (see Lemma \ref{A}), we have
\[E[v(\theta^\lambda(\theta_p)| \vee_n \mathcal{F}_{\theta_n}]=E[v(\theta^\lambda(\theta_p)|  \mathcal{F}_{\theta^-}] \quad \mbox{ a.s. on } \quad  A.\]
It follows that
\begin{equation}\label{eq-Ftmoins}
E[v(\theta^\lambda(\theta_p)|  \mathcal{F}_{\theta^-}]= v(\theta^-) \quad \mbox{ a.s. on } \quad  \quad \{\theta^\lambda(\theta_p)\ge \theta\} \cap A.
\end{equation}

Since $(v(\tau), \tau \in \mathcal{T}_{[\theta_p, \theta^\lambda(\theta_p)]})$  is a martingale (see Lemma 2.7), we have 
$E[v(\theta^\lambda(\theta_p)))|\mathcal{F}_\theta] = v(\theta)$ a.s. on ${\theta^\lambda(\theta_p) \ge  \theta}.$ 
Hence, by taking the condidional expectation with respect to $\mathcal{F}_{\theta^-}$ , we derive that
\[E[v(\theta^\lambda(\theta_p))|F\theta−]=E[v(\theta)|F\theta−]<v(\theta^-) \quad \mbox{  a.s. on } \quad B\cap A\cap \{\theta^\lambda(\theta_p)\ge \theta\},\]
which, with equality (4.8), yields that $B  \cap A \cap \{\theta^\lambda(\theta_p)  \ge \theta\} = \emptyset\}$ a.s.

 It follows that for each $p$, $\theta_p \le \theta^\lambda(\theta_p)<\theta$ a.s. on 
$B\cap A$ and $\theta^\lambda(\theta_p)\uparrow \theta$ a.s. on $B\cap A$.
In other words, the sequence $(\theta^\lambda(\theta_p))$ announces $\theta$ on $B\cap A$. By a property of the $(1-\lambda)$-optimal stopping times 
(see Lemma 2.5), for each $\lambda \in [0, 1[$
and for each $p \in  \mathbb{N}$, we have 
\begin{equation}
\lambda v(\theta^\lambda(\theta_p)) \le \phi(\theta^\lambda(\theta_p))	\quad \mbox{a.s.}(4.9) 
\end{equation} 
Note also that the sequence $(\theta^\lambda(\theta_p))$  clearly belongs to $\mathcal{A}(\theta, B\cap A)$. By letting $p$ tend to
$\infty$ and by using the LL property of the value function (see Theorem 4.3), we derive that
\[\lambda v(\theta^-) = \lambda \lim_{p\to \infty} v(\theta^\lambda(\theta_p)) \le  \limsup_{p\to \infty} \phi(\theta^\lambda(\theta_p)) \quad \mbox{a.s. on}   C = B\cap A \]
By definition of $\overline \phi_C(\theta)$, we also have 
\[ \limsup_{p\to \infty} \phi(\theta^\lambda(\theta_p)) \le  \phi_C (\theta)	\quad \mbox{ a.s.	on} \quad 	C.\]
 Hence, the previous inequalities lead to $\lambda v(\theta^-) \le  \limsup_{p\to \infty} \phi(\theta^\lambda(\theta_p)) \le  \overline \phi_C(\theta)$ a.s. on
$C$ and this holds for each $\lambda < 1$. By letting $\lambda$ tend to 1, we obtain 
\[v(\theta^-) \le  \limsup_{\lambda\uparrow 1} \limsup_{p\to \infty}  \phi(\theta^\lambda(\theta_p)) \le  \overline \phi_C (\theta) \quad \mbox{	a.s.	on}	C.\]
Moreover, since the value function is LL and since $ v \ge \phi,$ one can easily show that
$v(\theta^-) \ge \overline \phi_C (\theta)$ a.s. Hence, 
\begin{equation}v(\theta^-) = \limsup_{\lambda\uparrow 1} \limsup_{p\to \infty} \phi(\theta^\lambda(\theta_p)) = \overline \phi_C (\theta)\quad \mbox{	a.s.	on	}C.	(4.10)
\end{equation}
Suppose now that $(\phi(\theta),\theta \in  T0)$ is also LL at $\theta$ on $C$. Then, we have 
$\phi(\theta^−) = \phi_C(\theta)$ a.s. on $C$. Hence, the previous equality can be written $v(\theta^-) = \phi(\theta^-)$ a.s. on $C$. The proof is thus complete.
\end{proof}
\begin{remark}
By (\ref{pp}), we have  
\[ \overline \phi_C (\theta)=\limsup_{\lambda \uparrow 1} \limsup_{p \to \infty} \phi({\theta^{\lambda} (\theta_p)})   \,\,\, \,\,\,  {\rm a.s. }\,\,\,{\rm on}\,\,\,\,\,\, C=B\cap A.\]
Note that the above theorem is used in \cite{KQ}.
\end{remark}

\appendix{}

\section{A measurability property}

\begin{lemma}\label{measurability}
Let $S$ $\in $ $\mathcal{T}_0$ and $(S_n)$ be a non decreasing sequence in $\mathcal{T}_0$ such that $S_n \leq S$ a.s. for all $n$. Let $A:= \{S_n<S,\,\, \mbox{for all}\,\,n\, \}.$
\begin{itemize}
\item
If $f$ is an $\vee_n\mathcal{F}_{S_n}$-measurable real random variable, then $f \mathbb{1} _A$ is $\mathcal{F}_{S^-}$-measurable.
\item
Suppose that $\lim_{n \to \infty} S_n = S$ a.s.\,\\
If $f$ is an $\mathcal{F}_{S^-}$-measurable real random variable, then $f \mathbb{1} _A$ is $\vee_n\mathcal{F}_{S_n}$-measurable.
\item
Suppose that $\lim_{n \to \infty} S_n = S$ a.s.\,\\
For each non negative random variable $g$, we have 
\begin{equation}\label{eqr}
E[g \,|\, \vee_n\mathcal{F}_{S_n}] \,\mathbb{1} _A = E[g \,|\, \mathcal{F}_{S^-}]\, \mathbb{1} _A \,\, \,\,\mbox{a.s.}
\end{equation}
\end{itemize}
\end{lemma} 
\begin{proof} 
 Let us prove the first assertion. For this, it is sufficient to prove that this property holds for $f:= \mathbb{1} _B$, where $B$ $\in$ $\vee_n\mathcal{F}_{S_n}$. Let $\mathcal{G}:= \{B \in \vee_n\mathcal{F}_{S_n}, \,\,\, B \cap A \in\mathcal{F}_{S^-}\}$. First, $\mathcal{G}$ is a $\sigma$-algebra. Note that $ \vee_n\mathcal{F}_{S_n}$ is the 
$\sigma$-algebra generated by $\cup_n\mathcal{F}_{S_n}$. Now, for each $n$, if $B$ $\in$ $\mathcal{F}_{S_n}$, then $B \cap \{S_n < S \}$ $\in$ $\mathcal{F}_{S^-}$. It follows that $\mathcal{G}$ is a $\sigma$-algebra which contains $\cup_n\mathcal{F}_{S_n}$, which yields that $\mathcal{G}= \vee_n\mathcal{F}_{S_n}$. Hence, the first assertion holds.

Let us show the second one. For this, it is sufficient to prove that this property holds for $f:= \mathbb{1} _B$, where $B$ $\in$ $\mathcal{F}_{S^-}$. Let $\mathcal{G}':= \{B \in \mathcal{F}_{S^-}, \,\,\, B \cap A \in \vee_n\mathcal{F}_{S_n}\}$. First, $\mathcal{G}'$ is a $\sigma$-algebra. Recall that $\mathcal{F}_{S^-}$ is the $\sigma$-algebra generated by the set 
$\mathcal{C}:= \{C \cap \{t < S\},\,\,\, C \in\mathcal{F}_t \,\,\, \mbox{and}\,\, t \in \mathbb{R} _+\}$. Now, by using the assumption $\lim_{n \to \infty} S_n = S$ a.s.\,, one can show that if $B \in \mathcal{C}$, then $B \cap A \in \vee_n\mathcal{F}_{S_n}$. It follows that $\mathcal{G}'$ is a $\sigma$-algebra which contains $\mathcal{C}$, which yields that $\mathcal{G}'=\mathcal{F}_{S^-}$. Hence, the second assertion holds.

It remains to show the third one. By the first assertion, $E[g \,|\, \vee_n\mathcal{F}_{S_n}] \mathbb{1} _A$ is 
$\mathcal{F}_{S^-}$-measurable. Also, for each $B$ $\in$ $\mathcal{F}_{S^-}$, since by the second assertion 
$A \cap B \in \vee_n\mathcal{F}_{S_n}$, we get
\[E[E[g \,|\, \vee_n\mathcal{F}_{S_n}] \mathbb{1} _A \mathbb{1} _B] = E[g \mathbb{1} _A \mathbb{1} _B].\]
Hence, equality (\ref{eqr}) follows.
\end{proof} 

\section{Case of a reward process}
In this section, we consider the particular case where the reward is given by a progressive process $(\phi_t)_{0\le t\le T}$. By using the results provided in this paper and naturally some fine results of the General Theory of processes, we derive the corresponding results in the case of processes.

Let $(\phi_t)_{0\le t\le T}$ be a progressive process. The associated  family $(\phi_\theta, \theta\in \mathcal{T}_0)$ is then admissible.  Suppose that $v(0) < \infty$. 

Since the supermartingale family $(v^+(\theta),\theta\in \mathcal{T}_0)$ is RCE,  there exists a 
RCLL process $(v^+_t)_{0\le t \le T}$  that \emph{aggregates} the family $(v^+(\theta),\theta\in \mathcal{T}_0)$  that is such that
\begin{equation}\label{E.aggregv+}
v^+(\theta)= v^+_\theta \; \,\mbox{a.s.\; for all }\theta\in \mathcal{T}_0\,.
\end{equation}
This follows from a classical result (see for instance Therem 3.13 in Karatzas and Shreeve (1994) and Proposition 4.1 in Kobylanski et al. (2011)).

Define the process $(v_t)$ by
\begin{equation}\label{E.v}
v_t:= \phi_t\vee v^+_t.
\end{equation}
By Proposition \ref{prop.vv+}, we clearly have 
\begin{proposition} 
Suppose that the reward is given by a progressive process $(\phi_t)$ such that the associated value function satisfies
$v(0) < \infty$. Then, the adapted process $(v_t)$ defined by (\ref{E.v}) aggregates the value function family $(v(S), S\in \mathcal{T}_0)$, that is for all $S\in \mathcal{T}_0$, $v(S)=v_S\quad {\rm a.s.}$
\end{proposition}

\begin{remark} We point out that, in this work, we have only made the assumption $v(0) <+ \infty$, which is, in the case of a reward process, weaker than the assumption $(\phi_t)$ of class $\mathcal{D}$, required in the previous literature. 
\end{remark}

Note that according to the terminology of Dellacherie and Meyer (1980), the process $(v_t)$ a \emph{strong supermartingale} that is, a supermartingale such that the family 
$(v (\theta),   \theta \in \mathcal{T}_0)$ is a supermartingale family.  By a fine result of Dellacherie and Meyer (1980) (see Theorem 4 p408), it follows that there exists a right limited and left limited version of $(v_t)$, which we still denote by $(v_t)$. 
Note that $v^+_t = v_{t^+}$, $0 \leq t \leq T$, a.s.\,

If the family $(\phi_\theta,\theta\in \mathcal{T}_0)$ is right USCE, thanks to Theorem \ref{thm.eli}, for each $\lambda \in ]0,1[$, $\theta^{\lambda} (S)$ is $(1- \lambda)$-optimal for $v_S$. 

We now show that under some additional assumptions, for each $\lambda \in ]0,1[$, $\theta^{\lambda} (S)$ can be written as a hitting time of processes.

\begin{definition}
A process $(\phi_t)_{0\le t\le T}$ is said to be \emph{right upper semicontinuous} if for almost every $\omega$, the function $t \mapsto \phi_t (\omega)$ is right upper semicontinuous, that is for each $t \in 
[0,T]$,
\[\phi_t (\omega) \geq \limsup_{s  \to t^+} \phi_s (\omega).\]
\end{definition}

\begin{remark}
Note that if $(\phi_t)$ is right upper semicontinuous, then the associated family 
$(\phi_\theta,\theta\in \mathcal{T}_0)$ is right USC (along stopping times). If, moreover, $(\phi_t)$ is of class $\mathcal{D}$, then $(\phi_\theta,\theta\in \mathcal{T}_0)$ is right USCE. 
Our assumptions ``$v(0) < + \infty$ and $(\phi_\theta,\theta\in \mathcal{T}_0)$ right USCE'' are thus weaker than the classical assumptions ``$(\phi_t)$  right upper semicontinuous and of class $\mathcal{D}$'', made in the previous literature.  
\end{remark}


 For each $\lambda \in ]0,1[$ let us define the following stopping time for each $\omega$ by
\[\tau^{\lambda} (S)(\omega):= \inf\{t\geq S(\omega) \,, \lambda v_t(\omega)\leq \phi_t(\omega)\}.\]

\begin{proposition}
Suppose that $(\phi_t)_{0\le t\le T}$ is a right upper semicontinuous progressive process of class $\mathcal{D}$. Let $S$ $\in$ $\mathcal{T}_0$. Then, for each $\lambda \in ]0,1[$, 
 \[\tau^{\lambda} (S)= \theta^{\lambda} (S)\,\,\,\,{\rm a.s.}\]
\end{proposition}
\begin{proof} To simplify the notation, let us denote $\theta^{\lambda} (S)$ by $\theta^{\lambda}$ and $\tau^{\lambda} (S)$ by $\tau^{\lambda}$.

By Lemma \ref{stepun}, for almost every $\omega$ $\in$ $\Omega$, 
\begin{equation*}
\lambda v_{\theta^{\lambda} (\omega) }(\omega) \leq  \phi_{\theta^{\lambda} (\omega) } (\omega) \,\,{\rm a.s. \,},
\end{equation*}
which implies that $\tau^{\lambda} (\omega)$ $\leq$ $\theta^{\lambda} (\omega)$ by definition of $\tau^{\lambda} (\omega)$.

Let us show the other inequality. Suppose first we have shown that for each $\lambda \in ]0,1[$, the stopping time $\tau^{\lambda}$ satisfies
\begin{equation}\label{lemprim}
\lambda v_{\tau^{\lambda}}\leq  \phi_{\tau^{\lambda} } \,\,{\rm a.s. }
\end{equation}
Hence, $\tau^{\lambda}$ $\in$ ${\mathbb{T}}_S= \{\,\theta \in \mathcal{T}_S\, , \,\lambda v(\theta) \leq\phi(\theta) \,\,\mbox {a.s.} \,\}$ which implies that 
$\theta^{\lambda}= {\rm ess} \inf \,\, {\mathbb{T}}_S $ $\leq$ $\tau^{\lambda} $ a.s. 
Consequently, the desired equality $\tau^{\lambda} $ $=$ $\theta^{\lambda}$ a.s. follows.

It remains to show inequality (\ref{lemprim}). The proof is done by fixing $\omega$ $\in$ $\Omega$ such that the function $t \mapsto \phi_t (\omega)$ is right upper semicontinuous and the function 
$t \mapsto v^+_t (\omega)$ is right continuous. There exists a non increasing sequence of reals 
$(t^n)$ (which depend of $\omega$) in $[S(\omega),T]$ such that $\displaystyle{\tau^\lambda(\omega)= \lim_{n\to \infty} \downarrow t^n}$ and such that for each $n$,
$\lambda v_{t^n} (\omega)\leq \phi_{t^n}(\omega)$. 

Note that if $v_{\tau^\lambda}(\omega) = \phi_{\tau^\lambda}(\omega) $, it is clear that 
$\lambda v_{\tau^\lambda}(\omega) \leq  v_{\tau^\lambda}(\omega) \leq  \phi_{\tau^\lambda}(\omega)$.\\
Second, if $ v_{\tau^\lambda}(\omega) >  \phi_{\tau^\lambda}(\omega) $, then $ v_{\tau^\lambda}(\omega)=  v^+_{\tau^\lambda}(\omega)$.\\
Since the function $t \mapsto v^+_t (\omega)$ is right continuous, and since $t^n\downarrow \tau^\lambda(\omega)$, one has  $v^+_{\tau^\lambda}(\omega)=\displaystyle{\lim_{n\to\infty}} v^+_{t^n}(\omega)$.
Hence, 
\begin{equation*}
\lambda  v^+_{\tau^\lambda}(\omega)= \lim_{n\to \infty}v^+_{t^n}(\omega) \leq \liminf_{n\to \infty}\phi_{t^n}(\omega)\,\,\,
\end{equation*}
where the last inequality follows from the fact that, for each $n$, 
$\lambda v_{t^n} (\omega)\leq \phi_{t^n}(\omega)$.
Now, the right upper semicontinuous property of the function $t \mapsto \phi_t (\omega)$ yields that $\limsup_{n\to \infty}\phi_{t^n}
(\omega) \leq \phi_{\tau^\lambda}(\omega)$.\,It follows that if $ v_{\tau^\lambda}(\omega) >  \phi_{\tau^\lambda}(\omega) $, then 
\[\lambda v_{\tau^\lambda}(\omega)= \lambda  v^+_{\tau^\lambda}(\omega) 
\leq \limsup_{n\to \infty}\phi_{t^n}(\omega) \leq \phi_{\tau^\lambda}(\omega),\]
which makes the proof ended.
\end{proof} 
Suppose now that the reward process $(\phi_t)$ satisfies the assumptions of the above theorem and that $(\phi_{\theta},\theta\in \mathcal{T}_0)$ is left USCE. Then, the existence theorem (see Theorem \ref{thm.TAO}) can be applied. Moreover, we have:
\begin{proposition}\label{T.TAO3} 
Let $(\phi_t)_{0\le t\le T}$ be a right upper semicontinuous progressive process of class $\mathcal{D}$. Suppose that the associated family $(\phi_{\theta},\theta\in \mathcal{T}_0)$ is left USCE. Then, for all $S\in \mathcal{T}_0$, the stopping time  
$\tau_*(S)$ defined by 
\[\tau_*(S)= \inf\{\, t\geq S\,,\, v_t=\phi_t\,\}\]
satisfies $\tau_*(S)$ $=$ $\theta_{*}(S)$ a.s. In particular, $\tau_*(S)$ is the minimal optimal stopping time for $v_S$. 
\end{proposition}

\begin{proof}  Again, let us denote $\theta_{*}(S)$ by $\theta_{*}$, 
 $\tau_*(S)$ by $\tau_*$, $\theta^{\lambda} (S)$ by $\theta^{\lambda}$ and $\tau^{\lambda} (S)$ by $\tau^{\lambda}$.
 
First, we clearly have that $\,{\displaystyle  \lim_{\lambda\uparrow 1} \uparrow \tau^\lambda \leq \tau_*}$ a.s. Also, 
$\,{\displaystyle  \lim_{\lambda\uparrow 1} \uparrow \tau^\lambda = \lim_{\lambda\uparrow 1} \uparrow \theta^\lambda =  \theta_{*} }\,\,{\rm  a.s.,},$ 
and hence, $  \theta_* \leq \tau_*$ a.s.

Furthermore, for almost every $\omega$, since $v_{\theta_{*}(\omega)}(\omega)= \phi_{\theta_{*}(\omega)}(\omega)$, it follows that $ \tau_*(\omega) \leq \theta_{*}(\omega)$ by definition of 
 $ \tau_*(\omega)$.
 Thus, we have proven that $ \theta_{*} =  \tau_*  $ a.s.
\end{proof} 

In the sequel, we suppose that $(\phi_t)$ is of class $\mathcal{D}$ and that 
the family $( \phi _{\theta},   \theta \in \mathcal{T}_0)$ is right USCE.\\
Since $(v_t)$ a strong supermartingale of class $\mathcal{D}$, by a fine result of Mertens (see for example the second assertion of Proposition 2.26 in El Karoui (1980)), there exists a unique uniformly integrable RCLL martingale $(M_t)$, a unique predictable right continuous non decreasing process $(A_t)$ with $A_0=0$ and $E[A_T] < \infty$ and a unique right continuous adapted non decreasing process $(C_t)$, which is purely discontinuous with $C_0=0$ and $E[C_T] < \infty$, such that
\begin{equation*}
v_t = M_t - A_{t} - C_{t-}, \,\,\, 0 \leq t \leq T  \,\,\, {\rm a.s.} 
\end{equation*}
We have $\Delta C_t= v_t - v_{t^+}$.\\
Now, for each predictable stopping time $S$, we have  $E[v_S\,|\,\mathcal{F}_{S^-}]$ $=$ $^pv_S$ a.s.\,, where $(^{p}v_t)$ is the predictable projection of 
$(v_t)$. Hence, 
\begin{equation} \label{Meyetrois}
\Delta A_S = v_{S^-}- \,^{p}v_S.
\end{equation}
Recall that the process $(A_t)$ admits the following unique decomposition:
\[A_t = A^c_t + A^d_t,\]
where $(A^c_t)$ is the continuous part of $(A_t)$ and $(A^d_t)$ is its purely discontinuous part.
\par
Theorem \ref{thm.vphi} leads to the following property:

 \begin{proposition}\label{trei}
Suppose that the reward is given by a progressive process $(\phi_t)$ of class $\mathcal{D}$, such that 
the family $( \phi _{\theta},   \theta \in \mathcal{T}_0)$ is right USCE. 
   Let $\theta \in \mathcal{T}_0$. For almost every $\omega$ such that $\theta (\omega) < T$, if $v_{\theta}(\omega) \neq \phi_{\theta}(\omega)$, then the non decreasing function $s \mapsto A_{s}(\omega)$ is locally constant on the right of $\theta(\omega)$, that is there exists $t^{'}(\omega)$ $>$ $\theta(\omega)$, such that $A_{t^{'}(\omega) }(\omega)= A_{\theta}(\omega)$.
 \end{proposition}
\begin{proof} 
It is sufficient to show that for almost every $\omega$ s.t. $\theta (\omega) < T$,  if  for each $t^{'}$ $>$ $\theta (\omega)$, $A_{t^{'}}(\omega) > A_{\theta}(\omega)$, then $v _{\theta}(\omega) = \phi_{\theta}(\omega)$.

 Let us introduce the following set:
 \[ A =\left\{\omega \in \Omega   \mbox{s.t.} \,\,\theta (\omega) < T \,  \mbox{and} \left( \forall  t^{'}\in ]\theta(\omega), T], \, A_{t^{'}}(\omega) > A_{\theta}(\omega) \mbox{and  } v_{\theta}(\omega) > \phi_{\theta}(\omega) \right) \right\}.\]
 Without loss of generality, we can suppose for each $\omega$, the function $t \mapsto A_{t}(\omega)$ is right continuous.
 Then, one can easily see that for each $p \in \mathbb{N}^*$,
 \[ A = \bigcap_{n \geq p} \left\{\, A_{(\theta + \frac{1} {n}) \wedge T} > A_{\theta} \mbox{and  } v_{\theta} > \phi_{\theta}  \right\} \cap \{\theta <T\},\]
 which implies that $A$ $\in$ $\mathcal{F}_{\theta + \frac{1} {p}}$. Hence,
 $ A \, \in \, \bigcap_{p \geq 1}\mathcal{F}_{(\theta + \frac{1} {p})\wedge T}\,=  \, \mathcal{F}_{\theta},$
 by the right continuity of the filtration $(\mathcal{F}_t)$.\\
Suppose now that $P(A) >0$. 
 The definition of $A$ clearly yields that for each $\theta^{'} \in \mathcal{T}_0$ with $\mathbb{1} _A \theta  \leq  \mathbb{1} _A \theta ^{'}$ a.s. and  $\mathbb{1} _A \theta  \not \equiv  \mathbb{1} _A \theta ^{'}$ a.s., we have 
$\mathbb{1} _A  A_{\theta^{'}} \geq \mathbb{1} _A A_{\theta} $ a.s. and $\mathbb{1} _A  A_{\theta^{'}} \not \equiv \mathbb{1} _A A_{\theta} $ a.s. This implies that $v$ is a strict supermartingale on the right at $\theta$ on $A$. Thanks to Theorem \ref{thm.vphi}, we get $v _{\theta}(\omega) = \phi_{\theta}(\omega)$ a.s.\, on $A$, which provides the expected contradiction. 
Hence, we have $P(A) =0$.
\end{proof}

The following lemma hold.
\begin{lemma}\label{ggg}
Let $(u_t)$ be a strong supermartingale of class $\mathcal{D}$ and such that the family $( u _{\theta},   \theta \in \mathcal{T}_0)$ is LCE. Then, the non decreasing predictable process $(A_t)$ of the Mertens decomposition of $(u_t)$ is continuous.
\end{lemma}  

\begin{remark}
This property is stated in Dellacherie and Meyer (1980) (Theorem 10 p214) in the case of a RCLL supermartingale but it still holds in the general case. 
\end{remark}
\begin{proof} By Proposition \ref{prop.vLCE}, since $(\phi_\theta, \theta\in T_0)$ is  USCE, the value function  $(v_\theta,\theta\in T_0)$ is LCE. Thanks to the previous lemma, the result follows. \end{proof} 
  
From this lemma, we derive the following property.

\begin{proposition}  Suppose the reward is given by a progressive process $(\phi_t)$ of class $\mathcal{D}$ and such that 
the family $(\phi _{\theta},   \theta \in \mathcal{T}_0)$ is USCE. The non decreasing predictable process $(A_t)$ of the Mertens decomposition of the value function $(v_t)$ is then continuous. 
\end{proposition}
\begin{proof} 
 Since $(\phi (\theta),   \theta \in \mathcal{T}_0)$ is an USCE, the value function $(v_t)$ is LCE. Thanks to the previous lemma, the result follows.
\end{proof} 
Recall that by a fine result of the General Theory of processes (see Appendix in El Karoui (1980)), for each adapted process $(\phi_t)$, there exists a predictable process $(\overline \phi_t)$ such that 
\[\overline \phi_t = \limsup_{s \to t,\, s<t} \phi_s, \, \, 0 \leq t \leq T, \,\, {\rm a.s.}\]
Note that in the case where $(\phi_t)$ is left limited, $\overline \phi_{t}$ $=$ 
$\phi_{t^-}$.

We now state the following property.
\begin{proposition}\label{prop.flatdeuxbis} 
  Suppose the reward 
$\phi$ is given by a right upper semicontinuous progressive process $(\phi_t)$ of class $\mathcal{D}$.\\
We then have:
\begin{itemize}
\item
$\Delta C_t = \Delta C_t \mathbb{1} _{\{v_{t} = \phi_{t}\}}= (v_t - v_{t^+}) \mathbb{1} _{\{v_{t} = \phi_{t}\}}$ a.s.
\item
For almost every $\omega$, the nondecreasing continuous function $t \mapsto A^c_t(\omega)$ is ``flat'' away from the set 
$\mathcal{H} (\omega) := \{t \in [0,T] \,, \,v_t (\omega) = \phi_t (\omega)\}$ i.e. 
$\int_0^T \mathbb{1} _{\{v_t > \phi_t \} } dA^c_t = 0$ a.s.\,
\item 
$\Delta A_t^d = \Delta A_t^d \mathbb{1} _{\{v_{t^-} = \overline \phi_{t}\}}$ a.s.
\item
If the filtration is left quasicontinuous and if $\theta$ is a predictable jump of $(v_t)$, then we have:
\[\Delta v_{\theta}=  \Delta A^d_{\theta}\mathbb{1} _{\{v_{{\theta}^-} = \overline \phi_{{\theta}}\}} = (v_{\theta}  -  \overline \phi_{{\theta}} )\mathbb{1} _{\{v_{{\theta}^-} = \overline \phi_{{\theta}}\}} \quad
\,\mbox{a.s.}\]
\end{itemize}
\end{proposition}
 \begin{remark}\label{KE} 
In the case of a continuous reward $(\phi_t)$, this proposition corresponds to Theorem D13 in Karatzas and Shreve (1998).
\par
Concerning the general case, the above property can also be derived from the results stated in El Karoui (1981).
\end{remark}
\begin{proof} 
 By Theorem \ref{thm.vphi}, we have that for each stopping time $\theta$ $\in$ $\mathcal{T}_0$, 
$v_{\theta} - v_{\theta^+} = (v_{\theta} - v_{\theta^+}) \mathbb{1} _{\{v_{\theta} = \phi_{\theta}\}}$. As $\Delta C_{\theta}= v_{\theta} - v_{\theta^+}$, the first point clearly holds.
\par
The proof of the second point is based on Proposition \ref{trei} and on some analytic arguments. Note that these analytic arguments are the same as those used in the proof of Theorem D13 in Karatzas and Shreve (1998).
Without loss of generality, we can suppose for each $\omega$, the maps $t\mapsto v_t(\omega)$, $t\mapsto \phi_t(\omega)$ are right continuous and $t\mapsto A^c_t(\omega)$ is continuous.\\
Let us denote by $\mathcal{J}(\omega)$ the set on which the nondecreasing function $t \mapsto A^c_t(\omega)$ is ``flat'':
\begin{equation*}\label{Jo}
\mathcal{J}(\omega) : = \{t \in ]0,T[ \, , \,\, \exists \varepsilon >0\,\,\, \mbox{with} \,\,\, A^c_{t -\varepsilon} (\omega)= A^c_{t +\varepsilon}(\omega)  \}
\end{equation*}
The set $\mathcal{J}(\omega)$ is clearly open and hence can be written as a countable union of disjoint intervals: 
${\displaystyle \mathcal{J}(\omega)  = \cup_i ]\alpha_i (\omega), \beta_i(\omega)[}$.
We consider 
\[\mathcal{\hat J}(\omega)  = \cup_i [\alpha_i(\omega), \beta_i(\omega)[ = \{t \in [0,T[ \, , \,\, \exists \varepsilon >0\,\,\, \mbox
{with} \,\,\, A^c_{t } (\omega)= A^c_{t +\varepsilon}(\omega)  \}.\]
The nondecreasing function $t \mapsto A^c_t(\omega)$ is now ``flat'' 
on $\mathcal{\hat J}(\omega)$  which means now that  
$\int_0^T \mathbb{1} _{\mathcal{\hat J}(\omega) } dA^c_t (\omega)= \sum_{i} (A^c_ {\beta_i(\omega) }- A^c_ {\alpha_i(\omega) }) =0.$ 
We next show that for almost every $\omega$,
$\mathcal{H}^c(\omega) \subset \mathcal{\hat J}(\omega),$
which clearly provides the desired result.\\
Let us denote by $\mathbb Q$ the set of rationals. By Proposition \ref{trei} applied to constant stopping times $\theta:= t$, where $t$ 
$\in$ $\mathbb{Q}\cap [0,T[$, it follows that for a.e. $\omega$,
\begin{equation}\label{f}
\{t \in \mathbb{Q}\cap [0,T[ \,\,
{\rm s.t.} \,\, v_t (\omega) > \phi_t (\omega)\} \subset \mathcal{\hat J}(\omega).
\end{equation}
Let us now show that the desired inclusion 
\[\mathcal{H}^c(\omega)=\{t \in  [0,T[ \,\,
{\rm s.t.} \,\, v_t (\omega) > \phi_t (\omega)\} \subset \mathcal{\hat J}(\omega)\]
holds for a.e. $\omega$.
Fix $\omega$ such that (\ref{f}) holds and fix $t$ $\in$ $\mathcal{H}^c(\omega)$. 
Since $v_t (\omega) > \phi_t (\omega)$ and since
the maps $t\mapsto v_t(\omega)$ and $t\mapsto \phi_t(\omega)$ are right continuous, there exists a non increasing sequence of rationals $t_n(\omega)$ $\in$ $\mathbb{Q}\cap [0,T[ $ such that 
${\displaystyle t= \lim_{n\to \infty}\downarrow  t_n(\omega)}$ with 
$v_{t_n(\omega)}(\omega) > \phi_{t_n(\omega)} (\omega)$ for each $n$. 
Using the above inclusion (\ref{f}), the equality $\mathcal{\hat J}(\omega)  = \cup_i [\alpha_i(\omega), \beta_i(\omega)[$ and the fact that ${\displaystyle t= \lim_{n\to \infty}\downarrow  t_n(\omega)}$, we derive that there exist $i$ and $n_0$ (which both depend on $\omega$) such that for each $n$ $\geq$ $n_0$, $t_n(\omega)$ $\in$ $[\alpha_i(\omega), \beta_i(\omega)[$. 
It follows that the limit $t$ $\in$ $[\alpha_i(\omega), \beta_i(\omega)[$, which gives that $t$ $\in$ $\mathcal{\hat J}(\omega)$.
Hence, the inclusion $\mathcal{H}^c(\omega) \subset \mathcal{\hat J}(\omega)$ is proven, which ends the proof of the second point. 
\par
Let us now show the third point. 
Let $\theta$ be a jump time of $(A^d_t)$. Since $(A^d_t)$ is predictable, $\theta$ is a predictable stopping time. Note now that 
$ \{\Delta A^d_{\theta} >0 \}= \{E [ v({\theta})  \, |\, \mathcal{F}_{\theta^-} ] < v({\theta^-})\}$. Hence, by Corollary 
\ref{predictable}, we have that $v_{\theta^-} = \overline \phi_C (\theta)$ a.s. on $C=\{\Delta A^d_{\theta} >0 \}$, where $\overline \phi_C (\theta)$ is defined by (\ref{env}). Now, one can easily show that $\overline \phi_{\theta}$ $=$ $\overline \phi_C({\theta})$ a.s. on $C$. The proof of the third point is thus complete.
\par
It remains to show the last point. Now, when the filtration is left quasicontinuous, the martingale 
$M$ only admits inaccessible jumps (as seen in Proposition \ref{gauchebis}). Consequently, if $\theta$ is a jump time for $(v_t)$ which is a predictable stopping time, then it corresponds to a jump of the nondecreasing predictable process $(A^d_t)$. We thus have 
\[\Delta v_{\theta}= v_{\theta}- v_{\theta^-} =  \Delta A^d_{\theta}\mathbb{1} _{\{v_{{\theta}^-} = \overline \phi_{\theta}\}} = (v_{\theta}  -  \overline \phi_{\theta} )\mathbb{1} _{\{v_{{\theta}^-} = \overline \phi_{\theta}\}},\] 
which makes the proof ended. 
\end{proof}

\begin{remark} 
Suppose now that the reward 
$\phi$ is given by a right upper semicontinuous and left USCE progressive process $(\phi_t)$ of class $\mathcal{D}$. One can prove that for each $S \in \mathcal{T}_0$, the maximal optimal stopping time $\check{\theta}( S)$ satisfies that for almost every $\omega$, 
\[\check{\theta}( S)(\omega) = \inf \{t\geq S(\omega) \, , \, v_t(\omega) \neq M_t(\omega)\}\wedge T.\]
 This corresponds to a well-known result of the Optimal Stopping Theory (see El Karoui (1981)). 
\end{remark}


\begin{thebibliography}{9}


\bibitem{DM1}
C. Dellacherie\ and\ P.-A. Meyer, {\it Probabilit\'es et potentiel}, \'Edition enti\`erement refondue, Hermann, Paris, 1975. MR0488194 (58 \#7757)

 

\bibitem{DM2}
C. Dellacherie\ and\ P.-A. Meyer, {\it Probabilit\'es et potentiel. Chapitres V \`a VIII}, revised edition, Actualit\'es Scientifiques et Industrielles, 1385, Hermann, Paris, 1980. MR0566768 (82b:60001)

\bibitem{EK}
N. El Karoui, Les aspects probabilistes du contr\^ole stochastique, in {\it Ninth Saint Flour Probability Summer School---1979 (Saint Flour, 1979)}, 73--238, Lecture Notes in Math., 876 Springer, Berlin. MR0637471 (83c:93062)


\bibitem{KS1}
I. Karatzas\ and\ S. E. Shreve, {\it Brownian motion and stochastic calculus}, Graduate Texts in Mathematics, 113, Springer, New York, 1988. MR0917065 (89c:60096)

\bibitem{KS2} 
I. Karatzas\ and\ S. E. Shreve, {\it Methods of mathematical finance}, Applications of Mathematics (New York), 39, Springer, New York, 1998. MR1640352 (2000e:91076)

\bibitem{KQ}
M. Kobylanski, M.-C. Quenez and E. Roger de Campagnolle, Dynkin games in a general framework. 
arXiv:1202.1930 

\bibitem{KQRd}
M. Kobylanski, M.-C. Quenez\ and\ E. Rouy-Mironescu, Optimal double stopping time problem, C. R. Math. Acad. Sci. Paris {\bf 348} (2010), no.~1-2, 65--69. MR2586746 (2010m:60145)

\bibitem{KQR}
M. Kobylanski, M.-C. Quenez\ and\ E. Rouy-Mironescu, Optimal multiple stopping time problem, Ann. Appl. Probab. {\bf 21} (2011), no.~4, 1365--1399. MR2857451 (2012h:60130)

\bibitem{Maingueneau}
M. A. Maingueneau, Temps d'arr\^et optimaux et th\'eorie g\'en\'erale, in {\it S\'eminaire de Probabilit\'es, XII (Univ. Strasbourg, Strasbourg, 1976/1977)}, 457--467, Lecture Notes in Math., 649 Springer, Berlin. MR0520020 (81j:60055)


\bibitem{Neveu} 
J. Neveu, {\it Discrete-parameter martingales}, translated from the French by T. P. Speed., revised edition, North-Holland, Amsterdam, 1975. MR0402915 (53 \#6729)

\bibitem{Peskir} 
G. Peskir\ and\ A. Shiryaev, {\it Optimal stopping and free-boundary problems}, Lectures in Mathematics ETH Z\"urich, Birkh\"auser, Basel, 2006. MR2256030 (2008d:60004)

\bibitem{Shiryaev} 
A. N. Shiryayev, {\it Optimal stopping rules}, translated from the Russian by A. B. Aries, Springer, New York, 1978. MR0468067 (57 \#7906)
\end{thebibliography}
\end{document}